\documentstyle[12pt]{article}

\newtheorem{prop}{Proposition}[section]
\newtheorem{dfn}[prop]{Definition}
\newtheorem{theo}[prop]{Theorem}

\newtheorem{rem}[prop]{Remark}
\newtheorem{coro}[prop]{Corollary}
\newtheorem{exam}[prop]{Example}

\newtheorem{lem}[prop]{Lemma}

\def\C{{\bf C}}
\def\R{{\bf R}}
\def\Z{{\bf Z}}
\def\Q{{\bf Q}}
\def\P{{\bf P}}

\def\a{ \alpha }

\def\d{ \delta }
\def\e{ \varepsilon }
\def\g{ \gamma}

\def\L{ \Lambda }

\def\s{ \sigma }
\def\S{ \Sigma }

\def\p{ \varphi }

\def\ra{ \rightarrow }

\title{Manin's conjecture for toric varieties}

\author{
Victor V. Batyrev\thanks{Supported by Deutsche Forschungsgemeinschaft} \\
\small  Universit\"at-GHS-Essen, Fachbereich  6,  Mathematik \\
\small  Universit\"atsstr. 3,  45141  Essen, Germany  \\
\small  e-mail: victor.batyrev@aixrs1.hrz.uni-essen.de \\
and \\
Yuri Tschinkel\\
\small Dept. of Mathematics, U.I.C.\\
\small Chicago, (IL) 60608,  U.S.A.  \\
\small e-mail: yuri@math.uic.edu
}

\begin{document}

\date{}

\maketitle

\thispagestyle{empty}

\begin{abstract}
We prove an asymptotic formula conjectured by Manin for the number
of $K$-rational points  of bounded height
with respect to the anticanonical line bundle
for arbitrary smooth projective
toric varieties over a number field $K$.
\end{abstract}

\newpage

\tableofcontents

\bigskip
\bigskip
\bigskip

\vskip 0,5cm
{\Large {\bf Introduction}}

\bigskip

Let $X$ be a projective  algebraic variety defined over a
number field $K$, $X(K)$ the set of
$K$-rational points of $X$.
We are interested in relations between the geometry of $X$
and  diophantine properties of  $X(K)$ in the situation when
$X(K)$ is infinite. The main object of our study
is the {\em height function} on $X(K)$ with respect to a
metrized line bundle ${\cal L}$.
A metrized line bundle is a pair ${\cal L}=(L, \| \cdot \|_v) $,
consisting of a line bundle $L$
equipped with  a family $\{ \| \cdot \|_v \}$
of $v$-adic metrics ($v$ runs over the
set ${\rm Val}(K)$ of all valuation of $K$), satisfying
certain conditions. This defines a height function
$H_{\cal L}: X(K)\ra \R_{> 0}$ by
$$
H_{\cal L}(x):=\prod_{v\in {\rm Val}(K)}\|f(x)\|_v^{-1}
$$
where $f$ is a $K$-rational local section of
the line bundle $L$ not vanishing in $x\in X(K)$.

Let ${\rm Pic}(X)$ be the Picard group of $X$ and
$\L_{\rm eff}\subset {\rm Pic}(X)$ the cone of effective divisors.
Assume that  the  class of  $L$ is contained in the
interior of the cone of effective divisors
$\L_{\rm eff}\subset {\rm Pic}(X)_{\R}$.
In this case, some positive tensor  power of $L$
defines a birational map of $X$
into some projective space. We denote by
$U_L\subset X$ the Zariski open subset
such that the
restriction of the above birational map to $U_L$
is an isomorphism on its image.

For any
Zariski open subset  $U\subset U_L$
consider the {\em height zeta-function} defined by
the following series  \cite{bat.man,franke-manin-tschinkel}:
\[
Z_{{\cal L},U}(s) = \sum_{x \in U(K)} H_{\cal L}(x)^{-s}. \]
Then $Z_{{\cal L},U}(s)$ converges absolutely and uniformly
for ${\rm Re}(s)\gg 0$.
A Tauberian theorem relates
the analytic properties  of the zeta-function with the
asymptotic behaviour of the
number $N(U,{\cal L},B)$ of $K$-rational points $x \in U(K)$
with $H_{\cal L}(x) \leq B$ as $B \rightarrow \infty$.

We define
$$
\a({\cal L},U) = \inf \{ a \in \R \,\mid \,
Z_{{\cal L},U}(s) \hskip 0,2cm {\rm converges
\,\, for } \,\, {\rm Re}(s)>a\}.
$$
We say that $U$ does not contain
{\em ${\cal L}$-accumulating subvarieties}
if for any non-empty Zariski open subset $U'\subset U$
we have $\a({\cal L},U)=\a({\cal L},U')$.

One expects  a good accordance between the
geometry of $X$ and diophantine properties of the set of
$K$-rational points of $X$  which are contained in the complement
$U \subset X$ to some proper closed subvarieties in $X$.
Denote by
${\cal K}^{-1}$  the  metrized anicanonical
line bundle.
The following conjecture is due to
Yu. I. Manin \cite{franke-manin-tschinkel}:
\bigskip

{\em Let $X$ be a smooth projective variety over a number field $K$
whose anticanonical line bundle is ample $($i.e., $X$ is a Fano
variety). Assume that the set $X(K)$ of $K$-rational
points is Zariski dense. Let $U\subset X$ be the largest
Zariski open subset
which doesn't contain ${\cal K}^{-1}$-accumulating subvarieties. Then
$$
N(U,{\cal K}^{-1} ,
B) = c(X,{\cal K}^{-1},K)B(\log B)^{k-1}( 1 + o(1))  \,\,\,
{\rm for}\,\, B\ra \infty
$$
where $k$ equals the rank of the Picard group
${\rm Pic}(X)$ over $K$ and
$c(X,{\cal K}^{-1},K)$ is some positive constant which depends
on $X$, $K$ and the choice of $v$-adic metrics on the anticanonical
line bundle.  }

\bigskip

The above conjecture was refined by E. Peyre \cite{peyre}
who defined Tamagawa numbers of Fano varieties and
proposed an interpretation of $c(X,{\cal K}^{-1},K)$
in terms of these  numbers.

The conjecture of Manin was proved  for generalized flag
varities, complete intersections of small degree
and for some blow ups of projective spaces
\cite{franke-manin-tschinkel,bat.man,peyre}.

In this paper we prove Manin's conjecture and compute the constant
$c(X,{\cal K}^{-1},K)$ for arbitrary
smooth projective equivariant
compactifications of  algebraic
tori over number fields, i.e., for toric varieties
${\bf P}_{\S}$ associated with a Galois invariant
finite polyhedral fan $\S$ \cite{vosk1}.
We restrict ourselves to  the case when $U$ is the dense torus
orbit. It is easy to show that $U$
doesn't contain ${\cal L}$-accumulating subvarieties for any
metrized line bundle ${\cal L}$. On the other hand,
it might happen that there is some larger Zariski open
subset $U' \subset {\bf P}_{\S}$ which contains $U$ and which does not
contain ${\cal L}$-accumulating subvarieties. It is possible
to show that the the asymptotic
formula for $N(U',{\cal K}^{-1} ,
B)$ does not depend on the choice of such a Zariski open $U'$, but
we decided  to postpone the proof of this fact.

One of our main ideas for the computation
of the height zeta-function on a toric variety ${\bf P}_{\S}$
is to introduce some
canonical simultaneous metrizations
on all line bundles and to obtain a pairing
$$
H_{\S}(x,{\bf s}):\,\, T(K)\times {\rm Pic}({\bf P}_{\S})_{\C}\ra \C
$$
between the set of rational points
$T(K)\subset {\bf P}_{\S}(K)$
in the Zariski open subset $T$ and the complexified
Picard group, extending the usual height pairing between
$T(K)$ and ${\rm Pic}({\bf P}_{\S})$.
This allows to  extend the one-parameter
zeta-function to a function $Z_{\S}({\bf s})$ defined on the
complexified Picard group
${\rm Pic}({\bf P}_{\S})_{\C}$ and holomorphic
when the ${\rm Re}({\bf s})$ is contained in the interior
of the cone
$[{\cal K}^{-1} ] +  \L_{\rm eff}$, where $\L_{\rm eff}\subset
{\rm Pic}({\bf P}_{\S})_{\R}$ is  the cone
of effective divisors  of ${\bf P}_{\S}$.

The second step is to use the multiplicative
group structure on the torus
$T\subset {\bf P}_{\S}$.  With our choice of metrics,
the height zeta-function becomes a function on the
adelic group $T({\bf A}_K)$ invariant under the
closed subgroup $T(K){\bf K}_T$, where
${\bf K}_T\subset T({\bf A}_K)$
is the maximal compact subgroup.
The key idea is to use the {\em Poisson formula} in order to
obtain an integral representation for
$Z_{\S}({\bf s})$.

Our third step is the study of analytic properties of
$Z_{\S}({\bf s})$ using the above integral formula and properties of
${\cal X}$-functions of convex cones.

In section 1 we introduce notations and  basic notions from the
theory of toric varieties over non-closed fields.

In section 2 we recall the definitions of Tamagawa numbers
of algebraic tori and Tamagawa numbers of
algebraic varieties with metrized anticanonical bundle.

In section 3  we define simultaneous metrizations of
all line bundles on toric varieties, introduce the
height zeta-function and give formulas for local Fourier
transforms of heights.

In section 4 we prove the Poisson formula which yields
an integral representation of the
height zeta-function.

In section 5 we formulate basic properties of
${\cal X}$-functions of convex finitely
generated polyhedral cones.

In section  6 we prepare the necessary analytic tools.

And finally, in section 7 we prove our main
theorem:

\medskip
{\em Let ${\bf P}_{\Sigma}$ be a smooth
projective compactification of an algebraic
torus $T$ over $K$. Let $k$ be
the rank of ${\rm Pic}({\bf P}_{\Sigma})$.
Then there is only a finite number $N(T,
{\cal K}^{-1}, B)$
of $K$-rational points  $x \in T(K)$ having the anticanonical
height $H_{{\cal K}^{-1}}(x) \leq B$.
Moreover,
\[ N(T, {\cal K}^{-1}, B) =
\frac{\Theta(\Sigma )}{(k-1)!}
\cdot B ( \log B)^{k-1}(1 + o(1)),\;\; B \ra \infty,  \]
with  the constant
$\Theta(\Sigma )=\alpha ( {\bf P}_{\Sigma})
\beta ({\bf P}_{\Sigma})\tau_{\cal K}({\bf P}_{\Sigma})$, where:

$1.$ $\alpha ( {\bf P}_{\Sigma})$ is a constant
defined by the geometry of  the  cone of effective divisors
$\Lambda_{\rm eff} \subset
{\rm Pic}({\bf P}_{\Sigma})_{\bf R}$;

$2.$ $\beta ({\bf P}_{\Sigma})$ is  the order of the
non-trivial part of the Brauer group of  ${\bf P}_{\Sigma}$;

$3.$ $\tau_{\cal K}({\bf P}_{\Sigma})$ is
the Tamagawa number associated with
the metrized canonical sheaf on
${\bf P}_{\Sigma}$ $($as it was defined
by E.Peyre in {\rm \cite{peyre}}$)$. }
\medskip

Our results provide first examples
for asymptotics on unirational varieties which
are not rational and on varieties without weak
approximation, in general. A new phenomenon
is the appearance of the non-trivial part of
the Brauer group $Br({\bf P}_{\S})/Br(K)$ in the
asymptotic constant. We don't need to assume that
the anticanonical class of ${\bf P}_{\S}$ is ample (i.e.,
${\bf P}_{\S}$ is a toric Fano variety).

This paper is a continuation of our paper \cite{BaTschi},
where we proved the conjecture of Manin about the distribution
of $K$-rational points of bounded height
for the case of projective compactifications of anisotropic tori.
An equvariant  compactification ${\bf P}_{\S}$ of an anisotropic torus
$T$ is much simpler in many aspects:
the cone of effective divisors $\L_{\rm eff}$ is
always simplicial, all $K$-rational
points of ${\bf P}_{\S}$  are contained in $T$, and the
group $T({\bf A}_K)/T(K)$ is compact (this last property
significantly  simplifies the Poisson formula).

\bigskip

{\bf Acknowledgements. }
We would like to thank our teacher
Yu.I. Manin for  many years of  encouragement
and support.
We benefited from conversations with
W. Hoffmann and  B. Z. Moroz.

The work of the second author was
supported by a Junior Fellowship from
the Harvard
Society of Fellows. The paper was completed
while he was visiting the Max-Planck-Arbeitsgruppe
Zahlentheorie in Berlin, GHS-Essen and IHES. He would like to thank
these institutions for hospitality
and perfect working conditions.

\section{Algebraic tori and toric varieties}

Let $X_K$ be an algebraic variety defined over a number
field $K$ and $E/K$ a finite extension of number fields.
We will denote the set of $E$-rational points of $X_K$
by $X(E)$ and by $X_E$ the $E$-variety obtained from
$X_K$ by  base change. We sometimes omit the subscript in
$X_E$ if the respective field of definition
is clear from the context.  Let
${\bf G}_{m,E}= {\rm Spec}( E[x,x^{-1}])$ be
the multiplicative group scheme over $E$.

\begin{dfn}
{\rm A linear algebraic group $T_K$ is
called  a {\em  $d$-dimen\-sio\-nal algebraic torus} if
there exists a finite extension  $E/K$ such that
$T_{E}$ is isomorphic to
$({\bf G}_{m,E})^d$.
The field $E$ is called the
{\em splitting field }  of  $T$.

For any field $E$ we denote by $\hat{T}_E =
{\rm Hom}\,( T, E^*)$
the group of regular $E$-rational
characters of $T$.
}
\label{opr.tori}
\end{dfn}

\begin{theo} {\rm \cite{grothendieck,ono1,vosk}}
There is a contravariant equivalence
between the category of algebraic tori
defined over a number field $K$ and the category of
torsion free
${\rm Gal}(E/K)$-modules of finite rank over
${\bf Z}$.
The functors are given by
$$
M \ra T = {\rm Spec}(K\lbrack M \rbrack); \;\;
T\ra \hat{T}_E.
$$
The above contravariant equivalence is
functorial under field extensions of
$K$.
\label{represent}
\end{theo}

Let ${\rm Val}(K)$ be the set of all
valuations of a global field $K$. Denote by $S_{\infty}$
the set of archimedian valuations of $K$.
For any $v \in {\rm Val}(K)$, we
denote by $K_v$ the completion of
$K$ with respect to $v$.
Let $E$ be a
finite Galois extension of $K$.
Let ${\cal V}$ be an extension of $v$ to $E$,
$E_{\cal V}$ the completion
of $E$ with respect to ${\cal V}$. Then
\[ {\rm Gal}(E_{\cal V}/ K_v ) \cong G_v
\subset G, \]
where $G_{v}$ is the decomposition subgroup of
$G$ and $ K_v \otimes_K E \cong \prod_{{\cal V} \mid v} E_{\cal V}.  $
Let $T$ be an algebraic torus over $K$
with the splitting field $E$.
Denote by
$T(K_v)=$ the $v$-adic
completion of $T(K)$ and by $T({\cal O}_v)\subset T(K_v)$
its  maximal compact subgroup.

\begin{dfn}
{\rm Denote by $T({\bf A}_K)$ the adele group of $T$.
Define
\[T^1({\bf  A}_K) = \{ {\bf t} \in T({\bf A}_K) \, : \,
\prod_{v \in {\rm Val}(K)}
\mid m(t_v) \mid_v = 1, \; {\rm for \; all}\; m \in \hat{T}_K \subset M  \}.
\]
Let
\[ {\bf K}_T = \prod_{v \in {\rm Val}(K)} T({\cal O}_v), \]
be the maximal compact subgroup of $T({\bf A}_K)$.
}
\end{dfn}

\begin{prop} {\rm \cite{ono1}}
The groups $T({\bf A}_K)$, $T^1({\bf  A}_K)$, $T(K)$,
${\bf K}_T$ have the following properties:

{\rm (i)} $T({\bf A}_K)/T^1({\bf  A}_K) \cong {\bf R}^t$, where $t$
is the rank
of $\hat{T}_K$;

{\rm (ii)} $T^1({\bf  A}_K)/T(K)$ is compact;

{\rm (iii)} $T^1({\bf  A}_K)/ T(K)\cdot {\bf K}_T $
is isomorphic to the direct product of a finite group
${\bf cl}(T_K)$ and
a connected compact abelian topological group which dimension
equals the rank $r'$ of the group
of ${\cal O}_K$-units in $T(K)$;

{\rm (iv)} $W(T) = {\bf K}_T \cap T(K)$ is  a finite
group of all torsion elements in $T(K)$.
\label{subgroups}
\end{prop}

\begin{dfn}
{\rm
We define the following
cohomological invariants of the algebraic
torus $T$:
$$
h(T)={\rm Card}[H^1(G,M)],
$$
$$
{\rm III}(T)={\rm Ker}\, \lbrack
H^1(G, T(K)) \rightarrow \prod_{v\in {\rm Val}(K)}
 H^1(G_v, T(K_v)) \rbrack,
$$
$$
i(T)= {\rm Card}[{\rm III}(T)].
$$
\label{coh.inv}
}
\end{dfn}

\begin{dfn}
{\rm
Let $\overline{T(K)}$ be the closure
of $T(K)$ in $T({\bf A}_K)$ in the
{\em direct product topology}.
Define the {\em
obstruction group to weak approximation} as
$$
A(T)= T({\bf A}_K)/\overline{T(K)}.
$$
\label{weak0}
}
\end{dfn}

\begin{rem}
{\rm
It is known that over the splitting field $E$ one has
$A(T_E)=0$.
}
\end{rem}

Let us recall standard facts about toric varieties over
arbitrary fields
\cite{danilov,demasur,fulton,oda,BaTschi}.

\begin{dfn}
{\rm A finite set $\Sigma$ consisting of  convex rational polyhedral
cones in $N_{\bf R} = N \otimes {\bf R}$ is called a {\em
$d$-dimensional fan} if the following conditions are satisfied:

(i) every cone $\sigma \in \Sigma$
contains $0 \in N_{\bf R}$;

(ii) every face $\sigma'$ of
a cone $\sigma \in \Sigma$ belongs to $\Sigma$;

(iii) the intersection of any
two cones in $\Sigma$ is a face of
both cones. }
\end{dfn}

\begin{dfn}
{\rm A  $d$-dimensional fan $\S$ is called
{\em complete and regular} if the
following additional conditions are satisfied:

(i) $N_{\bf R}$ is the union of cones from $\Sigma$;

(ii)  every
cone $\sigma \in \Sigma$ is generated by a
part of a ${\bf Z}$-basis of
$N$.\\
We denote by $\Sigma(j)$ the set of all
$j$-dimensional cones in
$\Sigma$. For each cone $\sigma \in \Sigma$
we denote by
$N_{{\sigma}, \bf R}$ the minimal
linear subspace containing $\sigma$. }
\label{def.fan}
\end{dfn}

\noindent
Let $T_K$ be a $d$-dimensional  algebraic torus
over $K$
with splitting field $E$ and
$G = {\rm Gal}\, (E/K)$. Denote
by $M$  the lattice $\hat{T}_E$ and by
$N ={\rm Hom}\, (M, {\bf Z})$ the dual abelian group.

\begin{theo}
A complete regular $d$-dimensional fan $\S$ defines
a smooth
equivariant compactification ${\bf P}_{\S,E}$
of  the $E$-split algebraic torus $T_E$. The {\em
toric variety}  ${\bf P}_{\S,E}$ has the following
properties:

(i) There is a
$T_E$-invariant  open covering by affine subsets $U_{\sigma,E}$:
\[ {\bf P}_{\Sigma,E} = \bigcup_{ \sigma \in \Sigma} U_{\sigma,E}. \]
The affine subsets  are defined as
$U_{\sigma,E} = {\rm Spec}(E \lbrack M \cap \check{\sigma}
\rbrack$),  where $\check{\sigma}$
is the cone in $M_{\bf R}$ which is
dual to $\sigma$.

(ii) There is a representation of
${\bf P}_{\Sigma,E}$ as a disjoint
union of split  algebraic
tori $T_{\sigma,E}$ of dimension
${\rm dim}\, T_{\sigma,E} = d - {\rm dim}\, \sigma $:
\[ {\bf P}_{\Sigma,E} =  \bigcup_{ \sigma \in \Sigma } T_{\sigma,E}. \]
For each $j$-dimensional cone $\sigma \in \Sigma{(j)}$
we denote by
$T_{\sigma,E}$ the kernel of a
homomorphism $T_E \rightarrow
({\bf G}_{m,E})^j$ defined by a ${\bf Z}$-basis of
the sublattice $N \cap N_{{\sigma},{\bf R}} \subset N$.
\end{theo}

\noindent
To construct  compactifications of non-split tori
$T_K$ over $K$,
we need a complete fan $\Sigma$ of  cones
having an additional combinatorial structure: an {\em
action of the Galois group }
$G={\rm Gal}(E/K)$  \cite{vosk1}.
The lattice $M=\hat{T}_E$ is a $G$-module and
we have a representation $\rho: G\ra {\rm Aut}(M)$.
Denote by $\rho^*$ the induced dual
representation of $G$ in ${\rm Aut}(N)
\cong {\rm GL}(d,{\bf Z})$.

\begin{dfn}
{\rm  A complete  fan $\Sigma \subset
N_{\bf R}$ is called
{\em $G$-invariant} if for any $g
\in G$ and for any $\sigma \in \Sigma$, one
has  $\rho^*(g) (\sigma) \in \Sigma$.
Let $N^G$ (resp. $M^G$, $N_{\bf R}^G$,
$M_{\bf R}^G$ and $\S^G$) be the subset
of $G$-invariant elements
in $N$ (resp. in $M$, $N\otimes {\bf R}$,
$M\otimes {\bf R}$ and $\S$).
Denote by $\Sigma_G \subset N_{\bf R}^G$
the fan consisting of all possible
intersections $\sigma \cap N_{\bf R}^G$
where $\sigma$ runs over all cones in $\Sigma$.
}
\label{opr.invar}
\end{dfn}

\noindent
The following theorem is due to Voskresenski\^i \cite{vosk1}:

\begin{theo}
Let $\Sigma$ be a complete regular $G$-invariant
fan in $N_{\bf R}$. Assume that the complete toric variety
${\bf P}_{\Sigma,E}$ defined over the splitting field
$E$ by the $G$-invariant fan $\Sigma$ is projective.
Then there exists a unique complete
algebraic variety  ${\bf P}_{\Sigma,K}$ over $K$
such that its base extension
${\bf P}_{\Sigma,K} \otimes_{{\rm Spec} (K)}
{\rm Spec}(E)$ is isomorphic to
the toric variety ${\bf P}_{\Sigma,E}$.
The above isomorphism respects
the natural $G$-actions on ${\bf P}_{\Sigma,K} \otimes_{{\rm Spec}(K)}
{\rm Spec}(E)$ and ${\bf P}_{\Sigma,E}$.
\end{theo}

\begin{rem}
{\rm
Our definition of heights
and the proof of the analytic properties of height zeta
functions do not use the projectivity of respective
toric varieties.  We note that there exist non-projective compactifications of
split algebraic tori. We omit
the technical question of existence of non-projective
compactifications of non-split tori.
}
\end{rem}

We proceed to describe the algebraic geometric structure
of the variety ${\bf P}_{\S,K}$ in terms of the fan with
Galois-action. Let ${\rm Pic}({\bf P}_{\S,K}) $ be the
Picard group and
$\Lambda_{\rm eff}$ the cone in
${\rm Pic}({\bf P}_{\S,K})$ generated by
classes of effective divisors.
Let ${\cal K}$ be the canonical line
bundle of ${\bf P}_{\S,K}$.

\begin{dfn}
{\rm A continuous function $\varphi\; : \;
N_{\bf R} \rightarrow {\bf R}$ is called {\em $\Sigma$-piecewise linear}
if the restriction of $\varphi$ to every
cone $\sigma \in \Sigma$ is a linear function. It is
called {\em integral} if $\varphi(N) \subset {\bf Z}$.
Denote  the group of $\S$-piecewise linear integral functions by $PL(\S)$.
}
\end{dfn}

We see that the $G$-action on $M$  (and $N$)
induces a $G$-action on the free
abelian group $PL(\S)$.
Denote by $e_1, \ldots, e_n$ the primitive
integral generators
of all $1$-dimensional cones in $\Sigma$.
A function $\varphi\in PL(\S)$ is determined by
its values on $e_i,\, (i=1,...,n)$.
Let $T_{i}$
be the $(d-1)$-dimensional torus orbit
corresponding to the cone
${\bf R}_{\geq 0}e_i \in \Sigma(1)$ and
$\overline{T}_i$ the Zariski closure of $T_i$ in
${\bf P}_{\Sigma,E}$.

\begin{prop}
Let ${\bf P}_{\S,K}$ be a smooth toric variety
over $K$ which is an equivariant compactification
of an algebraic torus $T_K$
with splitting field $E$  and  $\S$
the corresponding complete regular fan
with $G={\rm Gal}(E/K)$-action. Then:

(i) There is an exact sequence
$$
0\ra M^G\ra PL(\S)^G\ra {\rm Pic}({\bf P}_{\S,K})\ra
H^1(G,M) \ra 0.
$$

(ii) Let
$$
\S(1)=\S_1(1)\cup ...\cup \S_r(1)
$$
be the decomposition of $\S(1)$ into a union
of $G$-orbits.
The cone of effective divisors
$\L_{\rm eff}$ is generated by classes
of $G$-invariant divisors
$$
D_j = \sum_{{\bf R}_{\geq 0}e_i \in \S_j(1)}
\overline{T}_i \,\,\, (j=1,...,r).
$$

(iii)
The class of the anticanonical line
bundle ${\cal K}^{-1}\in {\rm Pic}({\bf P}_{\S,K})$
is the class of the $G$-invariant piecewise linear
function $\varphi_{\S}\in PL(\S)^G$ given by
$\varphi_{\S}(e_j)=1$ for all $j=1,...,n$.
\label{nonsplit.geom}
\end{prop}

\begin{theo} {\rm \cite{vosk,ct}} Let $T$ be an
algebraic torus over $K$ with splitting field $E$.
Let ${\bf P}_{\S,K}$ be a complete
smooth equivariant compactification of $T$.
There is an exact sequence:
\[ 0\ra A(T) \ra
Hom (H^1(G,{\rm Pic}({\bf P}_{\S,E})),\Q/\Z)
\ra {\rm III}(T)\ra
0.
\]
\label{weak}
\end{theo}

\begin{rem}
{\rm
The group $H^1(G,{\rm Pic}({\bf P}_{\S,E}))$
is canonically isomorphic to the non-trivial
part of the Brauer group
${\rm Br}({\bf P}_{\Sigma,K})/{\rm Br}(K)$, where
${\rm Br}({\bf P}_{\Sigma,K}) = H^2_{\rm et}({\bf P}_{\Sigma,K}, {\bf G}_m)$.
This group appears
as the obstruction group to the Hasse principle
and weak approximation in \cite{manin,ct}.
}
\end{rem}

\begin{coro}
{\rm
Let $\beta({\bf P}_{\Sigma})$  be the cardinality of
$H^1(G,{\rm Pic}({\bf P}_{\S,E}))$. Then
\[ {\rm Card} \lbrack A(T) \rbrack =
\frac{\beta({\bf P}_{\Sigma})}{i(T)}. \]
\label{weak1}
}
\end{coro}

\section{Tamagawa numbers}

In this section we recall the definitions of Tamagawa
numbers of tori following A. Weil \cite{weil}
and of algebraic varieties with a
metrized canonical line bundle following E. Peyre \cite{peyre}.
The constructions of Tamagawa numbers depend on
a choice of a finite set of valuations
$S\subset {\rm Val}(K)$ containing archimedian
valuations and places of bad reduction, but the Tamagawa numbers
themselves do not depend on $S$.

Let $X$ be a smooth  algebraic variety over $K$, $X(K_v)$ the
set of $K_v$-rational points of $X$.
Then a choice of local analytic coordinates $x_1, \ldots, x_d$ on $X(K_v)$
defines a homeomorphism $\phi\,: \, U \ra K_v^d$
in $v$-adic topology between an
open subset $U \subset X(K_v)$ and  $\phi(U) \subset K_v^d$.
Let $dx_1 \cdots dx_d$ be the Haar measure on $K_v^d$ normalized
by the condition
\[ \int_{{\cal O}_v^d} dx_1 \cdots dx_d = \frac{1}{(\sqrt{\delta_v})^d} \]
where $\delta_v$ is the absolute different of $K_v$.
Denote by $dx_1 \wedge \cdots \wedge dx_d$ the standard
differential form on $K_v^{d}$. Then
$f = \phi^*(dx_1 \wedge \cdots \wedge dx_d)$ is a local analytic section of
the  canonical sheaf ${\cal K}$. If $\| \cdot \|$ is a $v$-adic
metric on  ${\cal K}$, then we obtain the $v$-adic measure on $U$ by
the formula
\[ \int_{U'} \omega_{{\cal K},v} =
\int_{\phi(U')}  \| f(\phi^{-1}(x))
\|_v dx_1 \cdots dx_d,
\]
where $U'$ is arbitrary open subset in $U$.
The  measure  $\omega_{{\cal K},v}$ does not
depend on the
choice of local coordinates
and extends  to a global
measure on $X(K_v)$ \cite{peyre}.

\begin{dfn}
{\rm
\cite{ono1} Let $T$ be an algebraic torus
defined over a number field $K$ with
splitting field $E$. Denote by
\[ L_S(s, T;E/K) =
\prod_{v \in {\rm Val}(K)} L_v(s, T ;E/K) \]
the Artin $L$-function corresponding to
the representation
\[  \rho \; :\; G= {\rm Gal}(E/K)
\rightarrow {\rm Aut}(\hat{T}_E)
 \]
and a finite set $S \subset {\rm Val}(K)$
containing all
archimedian valuations
and  all non-archimedian valuations
of $K$ which are
ramified in $E$.
By definition,
$L_v(s,T;E/K) \equiv 1$ if $v \in S$, $L_v(s,T;E/K)=
{\rm det}(Id - q^{-s}_v F_v)^{-1}$ if $v \not\in S$, where $F_v \in
{\rm Aut}(\hat{T}_E)$ is a representative of the Frobenius automorphism.
 }
\end{dfn}

Let $T$ be an algebraic torus of dimension $d$ and
$\Omega$ a $T$-invariant algebraic $K$-rational differential
$d$-form. The form $\Omega$ defines an isomorphism
of the canonical sheaf on $T$ with the structure sheaf on
$T$. Since the structure sheaf has a canonical metrization,
using the above construction, we obtain  a $v$-adic  measure
$\omega_{\Omega,v}$  on $T(K_v)$. Moreover, according to A. Weil \cite{weil},
we have
\[
 \int_{T({\cal O}_v)} \omega_{\Omega,v} =
\frac{{\rm Card} \lbrack T(k_v) \rbrack}{q^d_v} = L_v(1, T; E/K)^{-1}  \]
for all $v\not\in S$.
We put  $d\mu_v = L_v(1, T; E/K) \omega_{\Omega,v}$
for all $v\in {\rm Val}(K)$.
Then the  local measures $d\mu_v$ satisfy
$$
\int_{T({\cal O}_v)} d\mu_v = 1
$$
for all $v\not\in S$.

\begin{dfn}
{\rm
We define the {\em canonical measure}
on the adele group $T({\bf A}_K)$
$$
 \omega_{\Omega,S} = \prod_{v \in {\rm Val}(K)}
 L_v(1, T; E/K) \omega_{\Omega,v} =
 \prod_{v \in {\rm Val}(K)} d\mu_v.
$$
\label{can.meas}
}
\end{dfn}

\noindent
By the product formula,  $\omega_{\Omega,S} $
does not depend on the choice of $\Omega$.
Let ${\bf dx}$ be the standard Lebesgue measure
on $T({\bf
A}_K)/T^1({\bf  A}_K)$. There exists
a unique Haar measure
$\omega^1_{\Omega,S}$ on $T^1({\bf A}_K)$ such that $\omega^1_{\Omega,S}
{\bf dx} =
\omega_{\Omega,S}$.
\bigskip

We proceed to define {\em Tamagawa measures}
on algebraic varieties following E. Peyre \cite{peyre}.
Let $X$ be a smooth projective
algebraic variety  over $K$ with a metrized
canonical sheaf  ${\cal K}$. We assume that $X$  satisfies
the conditions
$h^1(X, {\cal O}_X) = h^2(X, {\cal O}_X) = 0$.
Under these
assumptions, the N{\'e}ron-Severi group $NS(X)$  (or, equivalently,
the Picard group ${\rm Pic}(X)$  modulo torsion)
over the algebraic closure
$\overline{K}$  is a discrete continuous ${\rm Gal}(\overline{K}/K)$-module
of finite rank over ${\bf Z}$.
Denote by  $T_{NS}$ the corresponding
torus under the duality from \ref{represent}
and by $E_{NS}$ a splitting field.

\begin{dfn}
{\rm \cite{peyre}
The {\em adelic Tamagawa measure}
$\omega_{{\cal K},S}$ on $X({\bf A}_K)$ is defined  by
$$
 \omega_{{\cal K},S} = \prod_{v \in {\rm Val}(K)}
L_v(1, T_{NS}; E_{NS}/K)^{-1}\omega_{{\cal K},v}. $$
}
\end{dfn}

\begin{dfn}
{\rm
Let $t$ be the rank of the group of
$K$-rational characters $\hat{T}_K$ of $T$.
Then the  {\em Tamagawa number of } $T$ is defined as
\[ \tau(T) = \frac{b_S(T)}{l_S(T)} \]
where
\[ b_S(T) = \int_{T^1({\bf  A}_K)/T(K)} \omega^1_{\Omega,S} , \]
\[ l_S(T) = \lim_{s \rightarrow 1} (s-1)^t L_S(s, T; E/K). \]}
\label{tamagawa1}
\end{dfn}

\begin{dfn}
{\rm \cite{peyre}
Let $k$ be the rank of the N{\'e}ron-Severi group of
$X$ over $K$, and $\overline{X(K)}$
the closure of $X(K) \subset X({\bf A}_K)$
in the
direct product topology.
Then the  {\em Tamagawa number} of $X$ is defined by
\[ \tau_{\cal K}(X) = \frac{b_S(X)}{l_S(X)} \]
where
\[ b_S(X) =
\int_{\overline{X(K)}} \omega_{{\cal K},S} \]
whenever the adelic integral converges, and
\[ l_S^{-1}(X) =
\lim_{s \ra 1} (s-1)^k L_S(s, T_{NS};
 E_{NS}/K). \]
}
\label{tamagawa2}
\end{dfn}

\begin{rem}
{\rm
Notice that there is a difference in the choice
of convergence factors for the
Tamagawa measure on
an algebraic variety $X$ and for
the Tamagawa measure on an algebraic torus
$T$. In the first case, we  choose $L_v^{-1}(1, T_{NS};E_{NS}/K)$
whereas in the second case one uses  $L_v(1, T; E/K)$.
This explains the difference in the
definitions of
$l_S(X)$ and $l_S(T)$.
}
\end{rem}

\begin{rem}
{\rm For a toric variety ${\bf P}_{\S}\supset T$
one can take $E_{NS}=E$, where $E$ is a splitting field
of $T$.
\label{ENS}
}
\end{rem}

\begin{rem}
{\rm
It is clear that in both definitions the Tamagawa numbers do not
depend on the choice of the finite subset
$S\subset {\rm Val}(K)$.
E. Peyre ($\cite{peyre}$) proves the existence
of the Tamagawa number for
Fano varieties by using
the Weil conjectures. The same method shows the existence
of the Tamagawa number for smooth complete varieties $X$ satisfying
the conditions $h^1(X, {\cal O}_X) = h^2(X, {\cal O}_X) = 0$.
}
\end{rem}

\begin{theo} {\rm \cite{ono2}} Let $T$ be an algebraic torus
defined over $K$.
The Tamagawa number $\tau (T)$ doesn't depend on
the choice of a splitting field $E/K$. We have
$$
\tau (T)=h(T)/i(T).
$$
The constants $h(T),i(T)$ were defined in \ref{coh.inv}.
\label{tamagawa}
\end{theo}

We see that the Tamagawa number of an algebraic torus
is a rational number. We have $\tau({\bf G}_m(K)) =1$.
The Tamagawa number of a Fano variety with a metrized
canonical line bundle is certainly not rational
in general. For $\P^1_{\Q}$ with our metrization
we have $\tau_{\cal K}(\P^1_{\Q}) =1/\zeta_{\Q}(2)$.

\begin{prop} {\rm \cite{BaTschi}}
One has
$$
\int_{\overline{T(K)}} \omega_{{\cal K},S} =
\int_{\overline{{\bf P}_{\S}(K)}} \omega_{{\cal K},S}.
$$
\label{two-integrals}
\end{prop}

\section{Heights and their Fourier transforms}

Let $\varphi \in PL(\Sigma)^G_{\bf C}$.
Using the  decomposition of $\S(1)$ into a union
of $G$-orbits
$$
\S(1)=\S_1(1)\cup ...\cup \S_r(1),
$$
we can identify $\varphi$ with a $T$-invariant divisor with
complex coefficients
$$
 D_{\varphi} = s_1 D_1 + \cdots + s_r D_r
$$
where $s_j = \varphi(e_j) \in {\C}$ and $e_j$ is a primitive lattice
generator of some cone $\s \in \S_j(1)$ $(j =1, \ldots, r)$.
It will be convenient to identify an element
$\varphi =\varphi_{\bf s}\in PL(\Sigma)^G_{\bf C}$ with the
vector ${\bf s} = (s_1, \ldots, s_r)$ of its complex
coordinates.

Let us  recall the definition of heights
on toric varieties from \cite{BaTschi}. For our purposes
it will be sufficient to describe the restrictions of heights to
the Zariski
open subset $T\subset {\bf P}_{\S,K}$.

\begin{prop} Let $v\in {\rm Val} (K) $ be a valuation
and $G_v\subset G$ the decomposition group of $v$.
There is an injective homomorphism
$$
\pi_v: T(K_v)/T({\cal O}_v)\hookrightarrow N_v,
$$
which is an isomorphism for all
but finitely many  $v\in {\rm Val}(K)$.
Here $N_v=N^{G_v}\subset N$ for
non-archimedian $v$ and $N_v=N_{\R}^{G_v}$ for
archimedian valuations $v$.
For every non-archimedian valuation we can identify
the image of $\pi_v$ with a sublattice of finite
index in $N_v$.
\label{pi-image}
\end{prop}

\begin{dfn}
{\rm Let ${\bf s } \in {\C}^r$ be a complex vector defining
a complex piecewise linear $G$-invariant function $\varphi
 \in PL(\Sigma)^G_{\bf C}$.
For any  point  $x_v \in T(K_v) \subset {\bf P}_{\Sigma}(K_v)$,
denote by $\overline{x}_v$ the image of
$x_v$ in $N_v$,  where $N_v$ is considered
as a canonical lattice  in the real
space $N_{\bf R}^{G_v}$ for non-archimedian
valuations  (resp. as the real Lie-algebra
$N_{\R,v}$ of $T(K_v)$ for archimedian valuations).
Define the {\em complexified
local Weil function}
$H_{\S,v}(x_v, {\bf s})$ by the formula
\[H_{\S,v}(x_v, {\bf s}) =
e^{\varphi(\overline{x}_v)\log q_v }\]
where $q_v$ is  the cardinality of the residue field
$k_v$ of $K_v$ if $v$ is non-archimedian
and $\log q_v = 1$ if
$v$ is archimedian.
}
\end{dfn}

\begin{theo} {\rm \cite{BaTschi}}
The complexified  local Weil function $H_{\S,v}(x_v, {\bf s}
)$ satisfies the
following properties:

{\rm (i)} $H_{\S,v}(x_v,{\bf s})$ is  $T({\cal O}_v)$-invariant.

{\rm (ii)} If ${\bf s} = 0$, then $H_{\S,v}(x_v,{\bf s}) = 1$
for all $x_v \in T(K_v)$.

{\rm (iii)} $H_{\S,v}(x_v, {\bf s} + {\bf s}') =
 H_{\S,v}(x_v,{\bf s}) H_{\S,v}(x_v,{\bf s}')$.

{\rm (iv)} If ${\bf s}=(s_1,...,s_r)\in {\bf Z}^r$,
then  $H_{\S,v}(x_v, {\bf s})$ is a classical local Weil
function  corresponding to
a  Cartier divisor $D_{\bf s} =
s_1 D_1 + \cdots + s_r D_r$ on ${\bf P}_{\S,K}$.
\label{local.f}
\end{theo}

\begin{dfn}
{\rm For a piecewise
linear function $\varphi_{\bf s} \in PL(\Sigma)^G_{\bf C}$
we define the {\em complexified
 height function on $T(K)\subset {\bf P}_{\S,K}(K)$} by
\[ H_{\Sigma}(x, {\bf s}) =
\prod_{v \in {\rm Val}(K)} H_{\S,v}(x_v, {\bf s}). \]}
\end{dfn}

\begin{rem}
{\rm
Although
the local heights are defined only
as functions on  $PL(\S)_{\C}^G \cong {\C}^r$, the
product formula implies that
for $x\in T(K)$ the global
complexified height function descends to
the Picard group ${\rm Pic}({\bf P}_{\S,K})_{\C}$.
Moreover, since $H_{\Sigma}(x, {\bf s})$
is the product of local complex
Weil functions $H_{\S,v}(x, {\bf s})$ and
since for all  $x_v \in T({\cal O}_v)$ we have
$H_{\S,v}(x_v, {\bf s}) = 1$ for  all
$v$,
we can immediately
extend $H_{\Sigma}(x,{\bf s})$ to a
function on
$T({\bf A}_K)\times PL(\S)^G_{\C}$.
}
\end{rem}

\begin{dfn}
{\rm Let
$ \S(1) = \S_1(1) \cup \cdots \cup \S_l(1) $
be the decomposition of $\S(1)$ into a disjoint union of $G_v$-orbits.
Denote by $d_j$ the length of the $G_v$-orbit $\S_j(1)$
$(d_1 + \cdots + d_l = n)$.
We establish a
1-to-1 correspondence $\S_j(1) \leftrightarrow u_j$ between
the $G_v$-orbits $\S_1(1), \ldots, \S_l(1)$ and independent variables
$u_1, \ldots, u_l$.
Let $\sigma \in \S^{G_v}$ be any
$G_v$-invariant cone and
$\S_{j_1}(1) \cup  \cdots \cup \S_{j_k}(1)$ the set  of all
$1$-dimensional faces of $\sigma$.
We define the rational function
$R_{\sigma}(u_1, \ldots, u_l)$ corresponding to $\sigma$ as follows:
\[ R_{\sigma}(u_1, \ldots, u_l) : =
\frac{u_{j_1}^{d_{j_1}} \cdots u_{j_k}^{d_{j_k}}
}{(1 - u_{j_1}^{d_{j_1}}) \cdots (1 - u_{j_k}^{d_{j_k}}) }. \]
Define the polynomial $Q_{\S}(u_1, \ldots, u_l)$ by the
formula
\[\sum_{\sigma \in \S^{G_v}} R_{\sigma}(u_1, \ldots, u_l) =
\frac{Q_{\S}(u_1, \ldots, u_l)}
{(1 - u_1^{d_1}) \cdots (1- u_l^{d_l}) }. \]
}
\end{dfn}

\begin{prop} {\rm \cite{BaTschi}}
Let $\Sigma$ be a complete  regular $G_v$-invariant
fan. Then the polynomial
\[ Q_{\S} (u_1, \ldots, u_l) - 1 \]
contains only monomials of degree $\geq 2$.
\label{p-function}
\end{prop}

Let $\chi$ be a topological character
of $T({\bf A}_K)$ such that
its $v$-component $\chi_v\, : \, T(K_v) \rightarrow S^1 \subset
{\bf C}^*$ is trivial on $T({\cal O}_v)$.
For each $ j \in \{ 1, \ldots, l\}$, we denote by
$n_j$ one  of $d_j$ generators of all  $1$-dimensional
cones of the $G_v$-orbit $\S_j(1)$; i.e., $G_vn_j$ is the set of
generators of $1$-dimensional cones in $\S_j(1)$.
Recall ($\ref{pi-image}$) that for
non-archimedian valuations,
$n_j$ represents an element
of $T(K_v)$ modulo $T({\cal O}_v)$. Therefore, $\chi_v(n_j)$ is
well-defined.
By ($\ref{pi-image}$) we know that the homomorphism
$$
\pi_v : T(K_v)/T({\cal O}_v)\ra N_v
$$
is an isomorphism for almost all $v$.
We call these valuations {\em good}.

\begin{dfn}
{\rm Denote by
$\hat{H}_{\Sigma,v} (\chi_v, -{\bf s})$ the value
at $\chi_v$ of the
{\it multiplicative}
Fourier transform of the local Weil function
$H_{\S,v}(x_v,-{\bf s})$ with
respect to the $v$-adic Haar measure
$d\mu_v$ on $T(K_v)$ normalized by
$\int_{T({\cal O}_v)} d\mu_v = 1$.}
\end{dfn}

\begin{prop} {\rm \cite{BaTschi}}
Let $v$ be a good non-archimedian
valuation of $K$ .
For any topological character
$\chi_v$ of $T(K_v)$ which is trivial on the
subgroup $T({\cal O}_v)$
and a piecewise linear function $\varphi = \varphi_{\bf s}
\in PL(\S)^G_{\C}$
one has
\[ \hat{H}_{\Sigma,v} (\chi_v, -{\bf s}) =
\int_{T(K_v)} H_{\S,v}(x_v,-{\bf s}) \chi_v(x_v)
 d\mu_v =  \]
 \[ =
\frac{Q_{\S}\left( \frac{\chi_{v}(n_1)}{{q_v}^{\varphi(n_1)}},
 \ldots, \frac{\chi_{v}(n_l)}{{q_v}^{\varphi(n_l)}} \right)}
{(1- \frac{\chi_{v}(n_1)}{{q_v}^{\varphi(n_1)}} )
\cdots (1 -
\frac{\chi_{v}(n_l)}{{q_v}^{\varphi(n_l)}} ) }.  \]
\label{integral.1}
\end{prop}

\begin{coro} {\rm \cite{BaTschi}}
{\rm
Let $v$ be a good non-archimedian valuation of $K$.
The  restriction of
 \[ \int_{T(K_v)} H_{\S,v}(x_v,-{\bf s})  d \mu_v \]
to the line $s_1 = \cdots = s_r = s$ is equal to
 \[ L_v( s,T; E/K)\cdot  L_v(s, T_{NS}, E/K) \cdot
 Q_{\S} (q_v^{-s}, \ldots, q_v^{-s}). \]
}
\label{loc-int}
\end{coro}

\begin{rem}
{\rm
It is difficult to calculate the
Fourier transforms of local heights for
the finitely many
"bad" non-archimedian
valuations $v$, because there is only
an embedding of finite index
$$
T(K_v)/T({\cal O}_v)\hookrightarrow N_v.
$$
However, for our purposes it will be sufficient
to use upper estimates for these local
Fourier transforms. One immediately sees that
for all non-archimedian valuations $v$
the local Fourier transforms
of $H_{\S,v}(x_v,-{\bf s})$ can be bounded
absolutely and uniformly in
all characters by a finite combination
of multidimensional geometric series in $q_v^{-1/2}$ in the domain
${\rm Re}({\bf s})\in \R_{>1/2}$.
\label{badreduction}
}
\end{rem}

Now we assume that $v$ is an archimedian  valuation.
By ($\ref{pi-image}$), we have
$T(K_v)/T({\cal O}_v) = N_{\bf R}^{G_v} \subset
N_{\bf R}$
where $G_v$ is the trivial group for the case
$K_v = {\bf C}$,
and $G_v = {\rm Gal}({\bf C}/ {\bf R}) \cong
 {\bf Z}/2{\bf Z}$
for the case $K_v = {\bf R}$.
Let $\langle \cdot,\cdot\rangle $ be the pairing
between $N_{\R}$ and $M_{\R}$ induced from the
duality between $N$ and $M$.
Let  $y$ be an arbitrary element of
the dual ${\bf R}$-space
$M_{\bf R}^{G_v} = Hom(T(K_v)/T({\cal O}_v), {\bf R})$.
Then  $\chi_y(x_v) = e^{- i \langle
 \overline{x}_v,y \rangle}$
is  a topological character
of $T(K_v)$ which is trivial on
$T({\cal O}_v)$.  We choose the Haar measure $d\mu_v$ on $T(K_v)$
as the product of the Haar measure $d\mu_v^0$
on $T({\cal O}_v)$ and
the Haar measure $d\overline{x}_v$ on
$T(K_v)/T({\cal O}_v)$.
We normalize the measures such that
the $d\mu_v^0$-volume of $T({\cal O}_v)$ equals $1$ and
$d\overline{x}_v$ is
the standard Lebesgue measure on $N_{\bf
R}^{G_v}$ normalized by the full sublattice $N^{G_v}$.

\begin{prop} {\rm \cite{BaTschi}}
Let $v$ be an archimedian
valuation of $K$.
The Fourier transform
$\hat{H}_{\S,v}(\chi_y,-{\bf s})$ of a
local archimedian Weil function
$$
H_{\S,v} (x_v,-{\bf s}) = e^{-\varphi_{\bf s}(\overline{x}_v)}
$$
is a rational function in
variables  $s_j = \varphi_{\bf s}(e_j)$ for
${\rm Re}({\bf s}) \in \R_{>0}$.
\label{archim.tr}
\end{prop}

{\it Proof.}
Let us consider the case  $K_v = {\bf C}$.
One has a decomposition of the space $N_{\bf R}$
into a union of $d$-dimensional cones $N_{\bf R} =
\bigcup_{\sigma \in \S(d)} \sigma$.
We calculate the Fourier transform as follows:
\[
\hat{H}_{\S,v}(\chi_y,-{\bf s}) =
\int_{N_{\bf R}} e^{-\varphi_{\bf s}(\overline{x}_v) -
 i \langle \overline{x}_v,y \rangle}
 d\overline{x}_v=
 \]
\[
= \sum_{\sigma \in \Sigma(d)} \int_{\sigma}
e^{-\varphi_{\bf s}(\overline{x}_v) -
 i \langle \overline{x}_v,y \rangle}
 d\overline{x}_v=
 \sum_{\sigma \in \Sigma(d)}
 \frac{1}{\prod_{e_j \in \sigma}
(s_j +   i \langle e_j,y \rangle)}.
\]
The case  $K_v = {\bf R}$ can be
reduced to the above situation.
\hfill $\Box$
\medskip

\section{Poisson formula}

Let ${\bf P}_{\S} $ be a toric variety and
$H_{\S}(x,{\bf s})$ the height function constructed
above.

\begin{dfn}
{\rm We define the zeta-function of the complex
height-function $H_{\Sigma}(x, {\bf s})$ as
\[ Z_{\Sigma}({\bf s}) = \sum_{x \in T(K)} H_{\Sigma}(x,{\bf -s}). \]}
\end{dfn}

\begin{theo}
The series $Z_{\Sigma}({\bf s})$
converges absolutely and uniformly for
 ${\bf s}$ contained in any compact in the
domain ${\rm Re}({\bf s})\in \R^r_{>1}$.
\label{convergence}
\end{theo}

\noindent
{\em Proof.}
It was proved in \cite{ono1} that we can always choose a finite set $S$
such that the natural map
$$
\pi_{S}\; : \; T(K) \ra \bigoplus_{v \not\in S} T(K_v)/T({\cal O}_v) =
\bigoplus_{v \not\in S} N_v
$$
is surjective. Denote by $T({\cal O}_S)$ the kernel of $\pi_S$ consisting
of all $S$-units in $T(K)$. Let $W(T) \subset T({\cal O}_S)$ the subgroup of
torsion elements in $T({\cal O}_S)$. Then $T({\cal O}_S)/W(T)$ has a natural
embedding into the finite-dimensional logarithmic space
$$
N_{\R,S} = \bigoplus_{v \in S} T(K_v)/T({\cal O}_v) \otimes {\R}
$$
as a sublattice of codimension $t$. Let $\Gamma$ be a full sublattice in
$N_{\R,S}$ containing the image of $T({\cal O}_S)/W(T)$. Denote by
$\Delta$ a bounded fundamental domain of $\Gamma$ in  $N_{\R,S}$.
For any $x \in T(K)$ we denote by $\overline{x}_S$ the image of $x$ in
$N_{\R,S}$. Define $\phi(x)$ to be the  element of $\Gamma$ such that
$\overline{x}_S - \phi(x) \in \Delta$. Thus, we have obtained
the mapping
$$
\phi\; : \; T(K) \ra \Gamma.
$$
Define a new function $\tilde{H}_{\S}(x, {\bf s})$ on $T(K)$ by
$$
\tilde{H}_{\S}(x, {\bf s}) = \prod_{v \in S} H_{\S,v}(\phi(x)_v,
{\bf s}) \prod_{v \not\in S} H_{\S,v}(x_v, {\bf s}).
$$
If ${\bf K}\subset \C^r$ is a compact in the domain
 ${\rm Re}({\bf s}) \in \R^r_{>1}$, then there exist two positive constants
 $C_1({\bf K}) < C_2({\bf K})$ such that
$$
0 < C_1({\bf K})
< \frac{\tilde{H}_{\S}(x, {\bf s})}{H_{\S}(x, {\bf s})} <
C_2({\bf K}) \;\; \mbox{\rm for ${\bf s} \in {\bf K}, \,x \in T(K)$},
$$
since $\overline{x}_v - \phi(x)_v$ belongs to some bounded subset $\Delta_v$
in $N_{\R,v}$ for any $x \in T(K)$ and $v \in S$.
Therefore, it is sufficient to prove that the series
$$
\tilde{Z}_{\S}({\bf s}) = \sum_{x \in T(K)} \tilde{H}_{\S}(x, - {\bf s})
$$
is absolutely converent for ${\bf s} \in {\bf K}$.
Notice that $\tilde{Z}_{\S}({\bf s})$ can be estimated from above by the
the following Euler product
$$
\left( \sum_{\g \in \Gamma} \prod_{v \in S}H_{\S,v}(\g_v, -{\bf s}) \right)
\prod_{v \not\in S} \left( \sum_{z \in N_v} H_{\S,v}(z, - {\bf s}) \right).
$$
The sum
$$
\sum_{\g \in \Gamma} \prod_{v \in S}H_{\S,v}(\g_v, -{\bf s})
$$
is an absolutely convergent geometric series for
${\rm Re}({\bf s}) \in \R^r_{>1}$.
On the other hand, the Euler product
$$
 \prod_{v \not\in S} \left( \sum_{z \in N_v} H_{\S,v}(z, - {\bf s}) \right)
$$
can be estimated from above by the product of zeta-functions
$$ C_3({\bf K}) \prod_{j =1}^r \zeta_{K_j}(s_j),$$
where $C_3({\bf K})$ is some constant depending on ${\bf K}$.
Since each $ \zeta_{K_j}(s_j)$ is absolutely convergent  for
${\rm Re}(s_j) > 1$, we obtain the statement.
\hfill
$\Box$

We need the Poisson formula in the following form:

\begin{theo}
Let ${\cal G}$ be a locally compact abelian group with
Haar measure $dg, {\cal H}\subset {\cal G} $ a closed
subgroup with Haar measure $dh$.
The factor group ${\cal G}/{\cal H}$ has a unique Haar measure $dx$
normalized by the condition $dg=dx\cdot dh$.
Let $F\,:\, {\cal G} \ra \R $ be an ${L}^1$-function on
${\cal G}$ and $\hat{F}$ its Fourier transform with respect
to $dg$. Suppose that $\hat{F}$ is also an ${L}^1$-function
on ${\cal H}^{\perp}$, where ${\cal H}^{\perp}$ is the group
of topological characters  $\chi \,: \, {\cal G} \ra S^1$
which are trivial on ${\cal H}$.
Then
$$
\int_{\cal H} F(x)dh=\int_{{\cal H}^{\perp}}\hat{F}(\chi) d\chi,
$$
where $d\chi$ is the orthogonal Haar measure on ${\cal H}^{\perp}$
with respect to the Haar measure $dx$ on ${\cal G}/{\cal H}$.
\label{poi}
\end{theo}

We will apply this theorem in the case when ${\cal G}=T({\bf A}_K)$
and ${\cal H}=T(K)$, $dg = \omega_{\Omega,S}$, and $dh$ is the
discrete measure on $T(K)$.

\begin{theo} (Poisson formula)
For all ${\bf s}$ with ${\rm Re}({\bf s})\in \R^r_{>1}$
we have the following formula:
$$
Z_{\Sigma}({\bf s})=\frac{1}{(2\pi )^t b_S(T)}
\int_{(T({\bf A}_K)/T(K))^*}
\left( \int_{T(A_F)}H_{\S}(x,-{\bf s})\chi(x)\omega_{\Omega,S}
\right) d\chi,
$$
where $\chi \in (T({\bf A}_K)/T(K))^* $
is a topological character of $T({\bf A}_K)$, trivial on
the closed subgroup $T(K)$ and
$d\chi$ is the orthogonal Haar measure on $(T({\bf A}_K)/T(K))^*$.
The integral converges
absolutely and uniformly to a holomorphic function in ${\bf s}$
in any compact in the domain ${\rm Re}({\bf s})\in \R^r_{>1}$.
\label{poiss}
\end{theo}

{\em Proof.} Because of \ref{convergence} we only
need to show that the Fourier transform
$\hat{H}_{\S}(\chi,-{\bf s})$ of the height function
is an $ L^1$-function on $(T({\bf A}_K)/T(K))^*$.
By \ref{integral.1} and uniform estimates at places of
bad reduction \ref{badreduction},
we know that the Euler product
$$
\prod_{v\not\in S_{\infty}}\hat{H}_{\S,v}(\chi_v,-{\bf s})
$$
converges absolutely and is
uniformly  bounded by a constant $c({\bf K})$ for all characters
$\chi$ and all ${\bf s}\in {\bf K}$, where ${\bf K}$ is  some
compact in the domain ${\rm Re}({\bf s})\in \R_{>1}^r$.

Since the height function
$H_{\S,v}(x,-{\bf s})$ is invariant under
$T({\cal O}_v)$ for all $v$, the  Fourier transform
of  $\hat{H}_{\S}(\chi,-{\bf s})$ equals
zero for characters $\chi$ which are non-trivial on the maximal
compact subgroup ${\bf K}_T$.
Denote by ${\cal P}$ the set of all such characters $\chi\in
(T({\bf A}_K)/T(K))^* $.

We have a non-canonical splitting
of  characters $\chi =\chi_l\cdot \chi_y$,
where $\chi_l\in (T^1({\bf A}_K)/T(K))^*$
and $\chi_y\in (T({\bf A}_K)/T^1({\bf A}_K))^*$.
Let us consider the logarithmic space
$$
N_{\R,\infty}=\bigoplus_{v\in S_{\infty}}T(K_v)/T({\cal O}_v)=
\bigoplus_{v\in S_{\infty}}N_{\R,v}.
$$
It contains the lattice $T({\cal O}_K)/W(T)$
of ${\cal O}_K$-integral points of
$T(K)$ modulo torsion.
Denote by
$
M_{\R,\infty}=\bigoplus_{v\in S_{\infty}}M_{\R,v}
$
the dual space.
It has a  decomposition as a direct sum of vector spaces
$M_{\R,\infty}=M_L\oplus M_Y$,
such that the space $M_L$ contains
the dual lattice $L:=(T({\cal O}_K)/W(T))^*$
as a full sublattice and the space $M_Y$ is isomorphic
to $(T({\bf A}_K)/T^1({\bf A}_K))^*= \hat{T}_K\otimes \R$.

By \ref{subgroups}, we have an exact sequence
$$
0\ra {\bf cl}^*(T)\ra {\cal P}\ra {\cal M}\ra 0,
$$
where ${\cal M}$ is the image of the projection of ${\cal P}$ to
$M_{\R,\infty}$ and ${\bf cl}^*(T)$ is a finite group.
We see that the character
$\chi\in {\cal P}$ is determined
by its archimedian component up to a finite choice.
Denote by $y(\chi)\in {\cal M}\subset M_{\R,\infty}$ the image
of $\chi\in {\cal P}$.

The following lemmas
will provide the necessary estimates of the
Fourier transform of local heights at archimedian places.
This allows to apply the Poisson formula \ref{poi}. \hfill $\Box $

\begin{lem} {\rm \cite{BaTschi}}
{\rm Let $\Sigma\subset N_{\R}$ be a
complete fan
in a real vector space of dimension $d$. Denote by $M_{\R}$ the
dual space. For all $m\in M_{\R}$
we have the following estimate
$$
|\sum_{\s \in \S (d)}
\frac{1}{\prod_{e_j\in \s}(s_j+i<e_j,m>)}|\le
\frac{1}{(1 + \|m\|)^{1+1/d}}.
$$
}
\end{lem}

\begin{coro}
{\rm
Consider
\[ \hat{H}_{\S,\infty}(\chi, -{\bf s}) =
\prod_{v \in S_{\infty}}
\hat{H}_{\S,v}(\chi, -{\bf s}) \]
as a function  on
\[ {\cal M}\subset M_{{\bf R}, \infty}
= \bigoplus_{v \in S_{\infty}} M_{{\bf R},v}. \]
Let $d'$ be the dimension of $M_{{\bf R}, \infty}$.
We have a direct sum decomposition
$M_{{\bf R}, \infty}=M_L\oplus M_Y$
of real vector spaces.
 Let $M'_Y\subset M_Y$ be any affine
subspace, $dy'$ the Lebesgue measure on $M'_Y$
and $L'\subset M_L$ any lattice.
Let  $g(y,-{\bf s})$ be a function on $M_{{\bf R}, \infty}$
satisfying the inequality $|g(y,-{\bf s})|\le c \|y\|^{\delta } $ for all
$y\in M_{{\bf R}, \infty}$, all ${\bf s}$ in some compact domain in
${\rm Re}({\bf s})\in \R_{>1/2}$, some $0<\delta < 1/d'$
and some constant $c>0$.  Then the series
\[
\sum_{y(\chi) \in L} \int_{M'_Y}
g(y,-{\bf s})\hat{H}_{\S,\infty}(y(\chi), -{\bf s}) dy'
\]
is absolutely and uniformly convergent to a holomorphic
function in ${\bf s}$ in this domain.
\label{infconver}
}
\end{coro}

{\em Proof.} We apply \ref{archim.tr} and observe that
on the space $N_{{\bf R}, \infty}$ we have a fan $\S_{\infty}$
obtained as the direct product of fans $\S^{G_v}$ for $v\in S_{\infty}$
(i.e., every cone in $\S_{\infty}$ is a direct product of cones in
$\S^{G_v}$).
\hfill $\Box $

\section{${\cal X}$-functions of convex cones}

Let $(A, A_{\bf R}, \L) $ be a triple consisting of
a free abelian group
$A$ of rank $k$, a $k$-dimensional real vector space
$A_{\bf R} = A \otimes {\bf R}$ containing $A$ as a sublattice of
maximal rank, and  a convex $k$-dimensional cone
$\Lambda \subset A_{\R}$ such that $\Lambda \cap - \L = 0
\in A_{\R}$. Denote by  $\L^{\circ}$ the interior  of $\L$ and
by  ${\L}_{\bf C}^{\circ} = {\L}^{\circ} + iA_{\R}$
the complex  tube domain over ${\L}^{\circ}$.
Let $( A^*, A^*_{\R}, \L^*) $ be the triple
consisting of the dual abelian group
$A^* = {\rm Hom}(A, \Z)$, the dual real vector space
$A^*_{\R} = {\rm Hom}(A_{\R}, \R)$, and the  dual cone
$\L^* \subset A^*_{\R}$.
We normalize the Haar measure $ {\bf d}{\bf y}$ on $A_{\R}^*$
by the condition:
${\rm vol}(A^*_{\R}/A^*)=1$.

\begin{dfn}
{\rm The {\em ${\cal X}$-function of}
${\L}$ is defined as
the integral
\[  {\cal X}_{\L}({\bf s}) =
\int_{{\L}^*} e^{- \langle {\bf s}, {\bf y}
 \rangle}  {\bf d}{\bf y}, \]
where ${\bf s} \in {\L}_{\bf C}^{\circ}$.  }
\end{dfn}

\begin{rem}
{\rm ${\cal X}$-functions of convex cones
have been investigated
in the theory of homogeneous cones
by M. K\"ocher, O.S. Rothaus, and
E.B. Vinberg \cite{koecher,rothaus,vinberg}. In these papers
${\cal X}$-functions were called {\em characteristic functions of
cones}, but we find such a notion rather misleading in view of the
fact that ${\cal X}_{\L}({\bf s})$ is the Fourier-Laplace transform
of the standard set-theoretic characteristic function of the
dual cone $\L^*$.  }
\end{rem}

\noindent
The  function  ${\cal X}_{\L}({\bf s})$
 has the following
properties \cite{rothaus,vinberg}:
\begin{prop}
{\rm (i)} If ${\cal A}$ is any invertible
linear operator on $\C^k$, then
\[ {\cal X}_{\L} ({\cal A}{\bf s})  = \frac{{\cal X}_{\L}({\bf s})}
{{\rm det}{\cal A}}; \]

{\rm (ii)} If ${\L} = {\bf R}^k_{\geq 0}$,  then
\[ {\cal X}_{\L}({\bf s}) = (s_1 \cdots s_k)^{-1},  \;{\rm for }
 \;{\rm Re}(s_i) > 0 ; \]

{\rm (iii)} If ${\bf s} \in {\L}^{\circ}$, then
\[ \lim_{{\bf s} \rightarrow \partial {\L}}
{\cal X}_{\L}({\bf s}) = \infty;  \]

{\rm (iv)} ${\cal X}_{\L}({\bf s}) \neq 0$ for all
${\bf s} \in {\L}_{\bf C}^{\circ}$.
\label{zeta.cone}
\end{prop}

\begin{prop}
If ${\L}$ is a $k$-dimensional
finitely generated polyhedral cone,
then ${\cal X}_{\L}({\bf s})$  is a rational function
$$
{\cal X}_{\L}({\bf s}) = \frac{P({\bf s})}{Q({\bf s})},
$$
where $P$ is a homogeneous polynomial,
 $Q$ is a product of all linear homogeneous forms defining
the codimension $1$ faces of
 $\L$, and ${\rm deg}\, P -
{\rm deg}\, Q = -k$.
\label{merom}
\end{prop}

\noindent
{\em Proof.} We subdivide the dual cone ${\L}^*$ into a finite
union of simplicial
subcones $\L_j^*$ $(j \in J)$.
Let $\L_j \subset A_{\R}$ be the dual cone
to $\L_j^*$. Then
$$
{\cal X}_{\L}({\bf s}) = \sum_{j \in J} {\cal X}_{\L_j}({\bf s}).
$$
By \ref{zeta.cone}(i) and (ii),
$${\cal X}_{\L_j}({\bf s}) = \frac{P_j({\bf s})}{Q_j({\bf s})} \;\;
( j \in J),$$
where $P_j$ is a homogeneous polynomial of degree $0$ and $Q_j$ is
the product of $k$ homogeneous linear forms defining the
codimension $1$ faces of $\L_j$.
Therefore, ${\cal X}_{\L}({\bf s})$ can be uniquely
represented up to constants
as a ratio of two homogeneous polynomials $P({\bf s})/Q({\bf s})$
with $g.c.d.(P,Q)=1$ where $Q$ is a product of linear homogeneous
forms defining some faces of $\L_j$ of codimension $1$. Since
$Q$ does not depend on a choice of a subdivision of $\L^*$ into a finite
union of simplicial cones $\L^*_j$, only linear homogeneous forms
which vanish on codimension $1$ faces of $\L$ can be factors of
$Q$. On the other hand, by \ref{zeta.cone}(iii), every linear homogeneous
form vanishing on a face of $\L$ of codimension $1$ must divide $Q$.
\hfill $\Box$

\begin{theo} Let $(A, A_{\R}, \L)$ and $(\tilde{A}, \tilde{A}_{\R},
\tilde{\L})$ be two triples as above, $k = {\rm rk}\, A$ and
$\tilde{k} = {\rm rk}\, \tilde{A}$, and  $\psi\;:\; A \ra \tilde{A}$
a homomorphism of free abelian groups with a finite cokernel
$A'$ (i.e., the corresponding
linear mapping of real vector spaces $\psi \;:\; A_{\R} \ra
\tilde{A}_{\R}$ is surjective), and $\psi(\L) = \tilde{\L}$.
Let $B= {\rm Ker}\, \psi \subset A$, ${\bf d}{\bf b}$ the Haar measure
on $B_{\R} = B \otimes {\R}$ normalized by the condition
${\rm vol}(B_{\R}/B)=1$.
Then for all ${\bf s}$ with
${\rm Re}({\bf s}) \in \Lambda^{\circ}$
the following formula holds:
$$
{\cal  X}_{\tilde{\L}}(\psi({\bf s}))
 = \frac{1}{(2\pi)^{k-\tilde{k}}|A'|}
\int_{B_{\R}} {\cal  X}_{{\L}}
({\bf s} + i {\bf b})  {\bf db},
$$
where $|A'|$ is the order of the finite abelian group $A'$.
\label{char0}
\end{theo}

{\em Proof.}  We have the dual injective homomorphisms
of free abelian groups
$\psi^*\;:\; \tilde{A}^* \ra A^*$ and of the corresponding
real vector spaces $\psi^*\;:\; \tilde{A}^*_{\R} \ra A^*_{\R}$.
Moreover, $\tilde{\L}^* = \L^* \cap \tilde{A}^*_{\R}$. Let
$C_{\L^*}({\bf y})$ be the set-theoretic characteristic function
of the cone $\L^* \subset A^*_{\R}$ and $C_{\L^*}(\tilde{\bf {y}})$ the
restriction of $C_{\L^*}({\bf y})$ to $\tilde{A}_{\R}^*$ which is
the set-theoretic characteristic function of $\tilde{\L}^* \subset
\tilde{A}_{\R}^*$. Then
$$
{\cal  X}_{\tilde{\L}}(\psi({\bf s})) =
\int_{\tilde{A}^*_{\R}} C_{\L^*}(\tilde{\bf {y}}) e^{- \langle
\psi({\bf s}), \tilde{\bf y}
 \rangle}  {\bf d}\tilde{\bf y}.
$$
Now we apply the Poisson formula to the last integral.
For this purpose  we notice that any additive topological character of
${A}^*_{\R}$ which vanishes on the subgroup $\tilde{A}^*_{\R} \subset
{A}^*_{\R}$ has the form
$$ e^{- i \langle {\bf b}, {\bf y} \rangle}, \;\;\; \mbox{ \rm where
 ${\bf b} \in B_{\R}$}.
$$
Moreover,
$$
\frac{{\bf db}}{(2\pi)^{k-\tilde{k}}|A'|}
$$
is the orthogonal Haar measure on $B_{\R}$ with respect to
the  Haar measures
${\bf d}\tilde{\bf y}$ and ${\bf d}{\bf y}$ on $\tilde{A}^*_{\R}$ and
$A^*_{\R}$ respectively. It remains to notice that
$${\cal  X}_{{\L}}
({\bf s} + i {\bf b}) =
\int_{{A}^*_{\R}} C_{\L^*}({\bf {y}}) e^{- \langle
{\bf s} + i{\bf b}, {\bf y}
 \rangle}  {\bf d}{\bf y}
 $$
is the value of the Fourier transform of
$C_{\L^*}({\bf {y}}) e^{- \langle {\bf s} , {\bf y} \rangle}$ on
the topological character of $A_{\R}^*/\tilde{A}_{\R}^*$ corresponding
to an element ${\bf b} \in B_{\R} \subset A_{\R}$.
\hfill $\Box$

\begin{coro}
{\rm
Assume that in the above situation ${\rm rk}\,= k - \tilde{k} =1$ and
$\tilde{A} = A/B$. Denote by $\g$ a generator of $B$. Then
$$
{\cal  X}_{\tilde{\L}}({\psi}({\bf s})) =
\frac{1}{2\pi i}\int_{{\rm Re}(z) = 0}
{\cal  X}_{\L}({\bf s} + z \cdot \g) dz,
$$
where $z = x + iy \in {\C}$.
\label{char1}
}
\end{coro}

\begin{coro}{\rm
Assume that a $\tilde{k}$-dimensional rational finite polyhedral cone
$\tilde{\Lambda} \subset \tilde{A}_{\R}$ contains exactly $r$ one-dimensional
faces with primitive lattice generators $a_1, \ldots, a_r \in \tilde{A}$.
We set $k := r$,  $A := {\Z}^r$ and denote by $\psi$
the natural homomorphism of lattices ${\bf Z}^r \ra \tilde{A}$
which sends the standard basis of ${\Z}^r$ into
$a_1, \ldots, a_r \in \tilde{A}$, so that $\tilde{\L}$ is the image
of the simplicial cone $\R^r_{\ge 0}\subset \R^r$
under the  surjective map of vector spaces $\psi\; : \; {\bf R}^r
\ra A_{\R}$.
Denote by $M_{\R}$ the kernel of $\psi$ and set $M := {\Z}^r \cap M_{\R}$.
Let ${\bf s}=(s_1,...,s_r)$ be the standard
coordinates in $\C^r$. Then
$$
{\cal X}_{\L}(\psi({\bf s}))=\frac{1}{(2\pi )^{r-k}|A'|}
\int_{M_{\R}}\frac{1}{\prod_{j=1,n}(s_j+iy_j)} {\bf d}{\bf y}
$$
where ${\bf dy}$ is the Haar measure on the additive
group $M_{\R}$ normalized
by the lattice $M$,
$y_j$ are the coordinates of ${\bf y}$ in $\R^r$, and
$|A'|$  is the index of the sublattice  in $\tilde{A}$ generated by
$a_1, \ldots, a_r$.
\label{int.formula}
}
\end{coro}

\begin{exam}
{\rm Consider an example of a non-simplicial
convex cone which
appears as the cone of effective divisors of the
split toric
Del Pezzo surface $X$ of anticanonical degree 6.
The cone ${\Lambda}_{{\rm eff}}$ has
6 generators corresponding to exceptional curves
of the first kind on $X$. We can construct
$X$ as the blow up of
3 points $p_1, p_2, p_3$ in general position
in ${\bf P}^2$.
Denote the  exceptional curves by $C_1, C_2, C_3,
C_{12}, C_{13}, C_{23}$, where $C_{ij}$ is
the proper pullback of the line
joining $p_i$ and $p_j$. Let
${\bf s} = s_1 [C_1] + s_2 [C_2] + s_3 [C_3]
 + s_{12}[C_{12}] +
s_{13}[C_{13}] + s_{23} [C_{23}] \in
\L_{\rm eff}^{\circ}$
be an element in the
interior of the cone of effective divisors.
The sublattice $M \subset {\bf Z}^6$ of rank $2$ consisting of
principal divisors is generated by
$\g_1 = C_1 + C_{13} - C_2 - C_{23}$ and $\g_2 = C_1 + C_{12} - C_3 -
C_{23} = 0$.

In our case, the integral formula in \ref{int.formula} is  a 2-dimensional
integral ($r =6$) which can be computed by applying twice the residue
theorem to two $1$-dimensional integrals like the one in \ref{char1}.
We  obtain the
following formula for the characteristic
function of ${\Lambda}_{\rm eff}$:  }
\[  {\cal X}_{\Lambda}(\psi({\bf s})) =
\frac{ s_1 + s_2 + s_3 + s_{12} + s_{13} + s_{23} }
{(s_1 + s_{23}) (s_2 + s_{13})(s_3 +
 s_{12})(s_1 + s_2 + s_3 )
(s_{12} + s_{13} + s_{23})}. \]
\end{exam}

\begin{dfn}
{\rm Let $X$ be a smooth proper algebraic variety.
Consider the triple $({\rm Pic}(X), {\rm Pic}(X) \otimes{\bf R},
\L_{\rm eff})$ where
$\L_{\rm eff} \subset {\rm Pic}(X)\otimes {\bf R}$  is
the cone generated by classes of effective
divisors on $X$.
Assume that the anticanonical
class $ \lbrack  {\cal K}^{-1}
\rbrack \in  {\rm Pic}(X)_{\bf R}$
is contained in the interior of $\L_{\rm eff}$. We define
the constant $\alpha(X)$ by
\[ \alpha(X) = {\cal X}_{\L_{\rm eff}}( \lbrack {\cal K}^{-1}
\rbrack). \]
}
\end{dfn}

\begin{coro}{\rm
If ${\L}_{\rm eff}$ is a finitely generated
polyhedral cone,
then $\alpha(X)$ is a rational number.
}
\end{coro}

\section{Some technical statements}

Let $E$ be a number field and
$\chi$  an unramified Hecke character
on ${\bf G}_m(A_E)$.
Its local components $\chi_v$ for all
valuations $v$ are given  by:
$$
\chi_v:
 {\bf G}_m(E_v)/{\bf G}_m({\cal O}_{v})\ra S^1
$$
$$
\chi_v(x_v)=|x_v|_v^{it_v}.
$$

\begin{dfn}
{\rm Let $\chi$ be an unramified Hecke character. We set
$$
y(\chi) : = \{ t_v \}_{v \in S_{\infty}(E)} \in
{\R}^{r_1 + r_2},
$$
where $r_1$ (resp. $ r_2$) is the number of real (resp. pairs of
complex) valuations of $E$. We also set
$$
\|  y(\chi) \| := \max_{v \in S_{\infty}(E)} |t_v|.
$$
\label{y-comp}
}
\end{dfn}
We will need uniform estimates for
Hecke $L$-functions in vertical strips.
They can be deduced using the Phragmen-Lindel\"of
principle \cite{rademacher}.

\begin{theo}
For any $\varepsilon > 0$
there exists a $\delta>0 $ such that for any
$0<\delta_1<\delta $ there exists a constant $c(\varepsilon,\delta_1) > 0$
such that the inequality
$$
| L_E(s, \chi) |
 \leq c(\varepsilon) ( 1 + |{\rm Im}(s)| + \|y(\chi)\| )^{\varepsilon}
$$
holds for all
$s$ with
$ \delta_1<|{\rm Re}(s) -1|< \delta$ and every Hecke L-function
$L_E(s, \chi)$ corresponding to an unramified
Hecke character $\chi$.
\label{estim}
\end{theo}

\begin{coro}{\rm
 For any $\varepsilon >0 $
there exists a $ \delta>0$ such
that for any compact
${\bf K}$ in the domain
$ 0<| {\rm Re}\,( s) - 1| <\delta$
there exists a
constant $C({\bf K},\varepsilon)$ depending only on
${\bf K}$ and $\epsilon$
such that
\[ | L_E(s,\chi) | \leq C({\bf K},\varepsilon)
(1 + \|y(\chi)\|)^{\varepsilon}
\]
for $s \in {\bf K}$ and every unramified
character $\chi$.
\label{m.estim}
}
\end{coro}

Let $\S$ be the Galois-invariant fan defining ${\bf P}_{\S}$ and
$\S(1)=\S_1(1)\cup ...\cup \S_r(1)$ the decomposition
of the set of one-dimensional generators of $\S$ into $G$-orbits.
Let $e_j$ be a primitive integral generator of $\sigma_j$,
$G_j \subset G$ the stabilizer of $e_j$. Denote by
$K_j \subset E$ the subfield of
$G_j$-fixed elements.
Consider the $n$-dimensional torus
$$
T'=\prod_{j=1}^r R_{K_j/K}({\bf G}_m).
$$
Let us recall the exact sequence
of Galois-modules from Proposition \ref{nonsplit.geom}:
$$
0\ra M^G\ra PL(\S)^G\ra {\rm Pic}({\bf P}_{\S})\ra H^1(G,M)\ra 0.
$$
It induces a map of tori $T'\ra T$ and a homomorphism
$$
a\,:\,\,T'({\bf A}_K)/T'(K)\ra T({\bf A}_K)/T(K).
$$
So we get a dual homomorphism for characters
$$
a^*:\,  (T({\bf A}_{K})/T(K))^*\ra \prod_{j =1}^r
({\bf G}_{m}({\bf A}_{K_j})/{\rm G}_m(K_j))^*.
$$

\begin{prop}{\rm \cite{drax1}}
The  kernel of $a^*$ is dual
to the obstruction group to weak approximation $A(T) $
defined in \ref{weak0}.
\label{kera}
\end{prop}

Let $\chi\in  (T({\bf A}_{K})/T(K))^*$ be a character. Then
$\chi\circ a$ defines $r$ Hecke characters of the idele groups
\[ \chi_j \; :\;
{\bf G}_m({\bf A}_{K_j}) \rightarrow S^1 \subset  {\bf C}^*, \;
j =1, \ldots, r. \]
If $\chi$ is trivial on ${\bf K}_T$, then all characters $\chi_j$
$(j =1, \ldots, r)$ are  trivial on the maximal
compact subgroups in ${\bf G}_m({\bf A}_{K_j})$.
We denote by  $L_{K_j}(s,\chi_j)$ the
Hecke $L$-function corresponding to the unramified
character $\chi_j$.

\begin{prop}
Let $\chi=(\chi_v)$ be a character and
$\hat{H}_{\S,v}(\chi_v,-{\bf s})$
the local Fourier transform
of the complex local height
function $H_{\S,v}(x_v,-{\bf s})$.
For any compact ${\bf K}$ contained
in the domain ${\rm Re}({\bf s})\in \R^r_{>1/2} $ there exists
a constant $c({\bf K})$ such that
\[ \prod_{v\not\in S}\hat{H}_{\S,v}(\chi_v,-{\bf s})\cdot
\prod_{i=1}^r L^{-1}_{K_j}(s_j,\chi_j) \le c({\bf K})\]
for all characters $\chi\in  (T({\bf A}_{K})/T(K))^*$.
\label{dmethod}
\end{prop}

The proof follows from explicit computations
of local Fourier transforms \ref{integral.1} and
is almost identical with the proof of
Proposition 3.1.3 in \cite{BaTschi}.

\begin{prop} There exists an $\varepsilon >0$ such that
for any open  $U\subset \C^r$
contained in the domain
$ 0< |{\rm Re}(s_j)-1|<\varepsilon $ for $j=1,...,r$
the integral
$$
\int_{(T({\bf A}_{K})/T(K))^*}\hat{H}_{\S}(\chi, -{\bf s})d\chi
$$
converges absolutely and uniformly to a holomorphic function
for ${\rm Re}({\bf s})\in U$.
\label{analytic}
\end{prop}
{\em Proof.} Using uniform estimates of Fourier transforms for
non-archimedian places of bad reduction (\ref{badreduction}) and
the proposition above we need only to consider the following
integral
$$
\int_{(T({\bf A}_{K})/T(K))^*}\hat{H}_{\S,\infty}(\chi,-{\bf s})
\prod_{j=1}^r L_{K_j}(s_j,\chi_j)d\chi.
$$
Observe that there exist constants $c_1>0$ and $c_2>0$
such that we have the following inequalities:
$$
c_1\|y(\chi)\|\le \sum_{j=1}^r\|y(\chi_j)\|\le c_2\|y(\chi)\|.
$$
Here we denoted by $\|y(\chi)\|$ the norm of $y(\chi)\in M_{\R,\infty}$.
Recall that since we only consider $\chi$ which are trivial on the
maximal compact subgroup ${\bf K}_T$, all characters $\chi_j$ are
unramified.
To conclude, we apply uniform  estimates
of Hecke L-functions from Corollary \ref{m.estim}
and the Corollary \ref{infconver}.
\hfill $\Box $

The rest of this section is devoted to
the proof of our main technical result.

Let ${\R}\lbrack {\bf s} \rbrack$
(resp. ${\C}\lbrack {\bf s} \rbrack$)
be the ring of polynomials in $s_1, \ldots, s_r$
with coefficients in
${\R}$ (resp. in ${\C}$),  ${\C}\lbrack \lbrack {\bf s}
\rbrack \rbrack$  the ring of formal power series in $s_1, \ldots, s_r$
with complex coefficients.

\begin{dfn}
{\rm Two elements $f({\bf s}),\, g({\bf s})\in {\C}\lbrack \lbrack {\bf s}
\rbrack \rbrack$ will be called {\em coprime}, if
$g.c.d.(f({\bf s}),\, g({\bf s})) =1$.
}
\end{dfn}

\begin{dfn}
{\rm Let  $f({\bf s})$ be an  element of ${\C}\lbrack \lbrack {\bf s}
\rbrack \rbrack$. By the {\em order} of a monomial
$s_1^{\a_1}...s_r^{\a_r}$ we mean the sum
of the exponents $\a_1+...+ \a_r$.
By  {\em multiplicity $\mu(f({\bf s}))$
of $f({\bf s})$ at
${\bf 0} = (0, \ldots, 0)$} we always mean
the minimal order of non-zero monomials appearing in the
Taylor
expansion of $f({\bf s})$ at ${\bf 0}$ .   }
\label{mult1}
\end{dfn}

\begin{dfn}
{\rm
Let  $f({\bf s})$ be  a  meromorphic at ${\bf 0}$ function.
Define the   {\em multiplicity $\mu(f({\bf s}))$ }
of $f({\bf s})$ at
${\bf 0}$ as
\[  \mu(f({\bf s})) = \mu(g_1({\bf s})) - \mu(g_2({\bf s})) \]
where $g_1({\bf s})$ and $g_2({\bf s})$ are two coprime
elements in  ${\C}\lbrack \lbrack {\bf s} \rbrack \rbrack$
such that $f = g_1/g_2$.
}
\label{mult2}
\end{dfn}

\begin{rem}
{\rm It is easy to show that for any two
meromorphic at ${\bf 0}$ functions  $f_1({\bf s})$ and $f_2({\bf s})$,
one has

(i) $\mu(f_1 \cdot f_2) = \mu(f_1) \cdot \mu(f_2)$ (in
particular, one can omit "coprime" in Definition \ref{mult2});

(ii) $\mu(f_1 + f_2) \geq \min \{ \mu(f_1),  \mu(f_2) \}$;

(iii)  $\mu(f_1 + f_2) = \mu(f_1)$ if $\mu(f_2) > \mu(f_1)$. }
\label{mult3}
\end{rem}

Using the  properties \ref{mult3}(i)-(ii), one immediately
obtains from  Definition \ref{mult1} the following:

\begin{prop}
Let $f_1({\bf s})$ and $f_2({\bf s})$ be two analytic at ${\bf
0}$ functions, $l({\bf s})$ a homogeneous linear function,
$\g = (\g_1, \ldots, \g_r) \in {\C}^r$ an arbitrary complex
vector with $l(\g) \neq 0$,
 and $g({\bf s}) := f_1({\bf s})/f_2({\bf s})$. Then
 the multiplicity  of
$$
 \tilde{g}({\bf s}): = \left(\frac{\partial}{\partial
z}\right)^k
g({\bf s} + z \cdot \g) |_{z = - l({\bf s})/l(\g)}
$$
at ${\bf 0}$ is at least $\mu(g) - k$.
\label{mult4}
\end{prop}

Let $\Gamma \subset {\Z}^r$ be a
sublattice, $\Gamma_{\R} \subset {\R}^r$ (resp.
$\Gamma_{\C} \subset {\C}^r$) the scalar extension of
$\Gamma$ to a {\R}-subspace
(resp. to a {\C}-subspace). We always assume that $\Gamma_{\R}
\cap {\R}_{\geq 0}^r = 0$. We set  $V_{\R}: = {\R}^r/\Gamma_{\R}$ and
$V_{\C}: = {\C}^r/\Gamma_{\C}$. Denote by $\psi$ the
canonical ${\C}$-linear projection ${\bf C}^r \ra V_{\C}$.

\begin{dfn}
{\rm A complex  analytic function $f({\bf s})=
f(s_1, \ldots, s_r): U \ra {\bf C}$ defined
on an open subset $U \subset {\bf C}^r$
is said to {\em descend to $V_{\C}$} if for any
vector $\a \in \Gamma_{\C}$ and any ${\bf u}= (u_1, \ldots, u_r) \in U$
one has
\[  f({\bf u}+ z \cdot \a) = f({\bf u}) \;\; \mbox{\rm for
all $\{ z \in \C\, :\,{\bf u}+ z \cdot \a \in U\}$}.\]
}
\end{dfn}

\begin{rem}
{\rm By definition, if $f({\bf s})$ descends to $V_{\C}$, then there exists
an analytic function $g$ on $\psi(U) \subset V_{\C}$ such that
$f = g \circ \psi$. Using Cauchy-Riemann equations, one
immediatelly   obtains that $f$ descends to $V_{\C}$ if and only
if for any vector $\a \in \Gamma_{\R}$ and
any ${\bf u}= (u_1, \ldots, u_r) \in U$, one has
\[ f({\bf u}+ iy \cdot \a)  = f({\bf u})\; \;
 \mbox{\rm for all $\{ y \in \R\, :\,{\bf u}+ iy \cdot \a \in
U\}$}. \]}
\label{desc}
\end{rem}

\begin{dfn}
{\rm An analytic function $W({\bf s})$ in the domain
${\rm Re}({\bf s}) \in {\R}_{>0}^r$ is called
{\em good with respect to $\Gamma$} if it satisfies the following conditions:

{(i)} $W({\bf s})$ descends to $V_{\C}$;

{(ii)} There exist pairwise coprime linear
homogeneous polynomials
$$
l_1({\bf s}), \ldots,
l_p({\bf s}) \in {\R}\lbrack {\bf s} \rbrack$$
and positive integers
$k_1, \ldots, k_p$ such that  for every $j \in \{1, \ldots, p \}$
 the
linear form $l_j({\bf s})$ descends to $V_{\C}$,  $l_j({\bf s})$
does not vanish for ${\bf s} \in {\R}_{>0}^n$, and
$$
P({\bf s}) = W({\bf s}) \cdot \prod_{j =1}^p l_j^{k_j}({\bf s})
$$
is analytic at ${\bf 0}$.

(iii) There exist a non-zero complex number $C(W)$ and
 a decomposition of $P({\bf s})$ into the sum
$$
P({\bf s}) = P_0({\bf s}) + P_1({\bf s})
$$
so that  $P_0({\bf s})$ is a homogeneous polynomial of degree $\mu(P)$,
$P_1({\bf s})$ is an  analytic function at ${\bf 0}$ with
$\mu(P_1) > \mu(P_0)$,   both
functions $P_0$, $P_1$ descend  to $V_{\C}$, and
$$
 \frac{P_0({\bf s})}{\prod_{j =1}^p l_j^{k_j}({\bf s})} =
C(W) \cdot {\cal  X}_{\L}(\psi({\bf s})),
$$
where  ${\cal  X}_{\L}$ is the ${\cal X}$-function of the cone $\L =
\psi({\R}^r_{\geq 0}) \subset V_{\C}$;
}
\end{dfn}

\begin{dfn}
{\rm If $W({\bf s})$ is a good with respect to $\Gamma$ as above, then the
meromorphic function
$$
\frac{P_0({\bf s})}{\prod_{j =1}^p l_j^{k_j}({\bf s})}
$$
will be called the {\em principal part of $W({\bf s})$ at ${\bf 0}$}
and the non-zero constant $C(W)$ the {\em principal coefficient
of $W({\bf s})$ at ${\bf 0}$}. }
\end{dfn}

Suppose  that $\Gamma \neq {\Z}^r$.
Let  $\g \in {\Z}^r$ be an element which is
not contained in $\Gamma$, $\tilde{\Gamma}: = \Gamma \oplus \Z < \g >$,
$\tilde{\Gamma}_{\R} := \Gamma_{\R} \oplus \R < \g >$, $\tilde{V}_{\R} :=
{\R}^r /\tilde{\Gamma}_{\R}$ and $\tilde{V}_{\C} :=
{\R}^r /\tilde{\Gamma}_{\C}$.

The following easy statement will be helpful in the sequel:

\begin{prop}
Let $f({\bf s})$ be an analytic at ${\bf 0}$ function,
$l({\bf s})$ a homogeneous linear function such that $l(\g) \neq
0$. Assume
that $f({\bf s})$ and $l({\bf s})$ descend to $V_{\C}$. Then
$$
\tilde{f}({\bf s}) : = f\left({\bf s} - \frac{l({\bf s})}{l(\g)}
\cdot \g\right)
$$
descends to $\tilde{V}_{\C}$.
\label{desc2}
\end{prop}

\begin{theo}
Let $W({\bf s})$ be a good function with respect to $\Gamma$  as above,
$$
\Phi({\bf s}) = \prod_{j\;:\; l_j(\g)=0} l_j^{k_j}({\bf s})
$$
the product of those linear
forms $l_j({\bf s})$ $j \in \{ 1, \ldots, p\}$ which vanish on $\g$.
Assume that $\tilde{\Gamma}_{\R} \cap {\R}_{\geq 0}^r = 0$
 and  the following statements hold:

{\rm (i)} The integral
$$
\tilde{W}({\bf s}) : = \int_{{\rm Re}(z) = 0} W({\bf s} + z \cdot \g) dz
, \;\; z \in \C
$$
converges  absolutely and uniformly
on any compact  in the domain ${\rm Re}({\bf s}) \in {\R}_{>0}^r$;

{\rm (ii)}
There exists $\d > 0$  such that the integral
$$
 \int_{{\rm Re}(z) = \d}
\Phi({\bf s}) \cdot W({\bf s} + z \cdot \g) dz
$$
converges absolutely and uniformly
in an  open neighborhood of ${\bf 0}$. Moreover, the
multiplicity of the meromorphic function
$$
\tilde{W}_{\d}({\bf s}): =   \int_{{\rm Re}(z) = \d}
 W({\bf s} + z \cdot \g) dz
$$
at  ${\bf 0}$ is at least $1 + {\rm rk}\, \tilde{\Gamma} - r$;

{\rm (iii)} For any ${\rm Re}({\bf s}) \in {\R}_{>0}^r$, the function
$$
\phi(t,{\bf s}) =
\sup_{0\leq {\rm Re}(z) \leq \d,\,{\rm Im}(z)=t }|W({\bf s}+ z\cdot \g)|
$$
tends to $0$ as $|t| \ra + \infty$.

Then $\tilde{W}({\bf s})$ is a good function with
respect to $\tilde{\Gamma}$, and $C(\tilde{W}) = 2 \pi i \cdot
C(W)$.
\label{desc3}
\end{theo}

\noindent
{\em Proof.} Assume that $l_j(\g) < 0$ for $j=1, \ldots, p_1$,
$l_j(\g) = 0$ for $j=p_1 +1, \ldots, p_2$, and
$l_j(\g) > 0$ for $j=p_2 +1, \ldots, p$. In particular, one has
\[ \Phi({\bf s}) = \prod_{j = p_1 + 1}^{p_2} l_j^{k_j}({\bf s}). \]
Denote by $z_j$ the solution of
the equation
\[
l_j({\bf s}) + z l_j(\g) = 0,\;\;j =1, \ldots, p_1.
\]

Let $U$ be the intersection of ${\R}^r_{>0}$ with an open
neighborhood of ${\bf 0}$ where
$$ \Phi({\bf s}) \cdot
\tilde{W}_{\d}({\bf s})
$$
is analytic. Then  both functions
$\tilde{W}_{\d}({\bf s})$ and $\tilde{W}({\bf s})$
are analytic in $U$. Moreover, the
integral formulas for $\tilde{W}_{\d}({\bf s})$ and $\tilde{W}({\bf s})$
show that the equalities
$\tilde{W}_{\d}({\bf u}+ iy \cdot \g) =\tilde{W}_{\d}({\bf u})$
and $\tilde{W}({\bf u}+ iy \cdot \g) =\tilde{W}({\bf u})$ hold
for any $y \in \R$ and ${\bf u},{\bf u}+ iy \cdot \g \in U$. Therefore,
both functions
$\tilde{W}_{\d}({\bf s})$ and $\tilde{W}({\bf s})$
descend to $\tilde{V}_{\bf C}$ (see Remark \ref{desc}).

Using assumptions (i)-(iii) of Theorem,  we
can apply the residue theorem and obtain
\[  \tilde{W}({\bf s})  - \tilde{W}_{\d}({\bf s}) = 2 \pi i\cdot
\sum_{j=1}^{p_1} {\rm Res}_{z = z_j} W({\bf s} + z \cdot \g)\]
for ${\bf s} \in U$.

We denote by $U(\g)$ the open subset of $U$ defined by the
inequalities
$$
\frac{l_j({\bf s})}{l_j(\g)} \neq \frac{l_{m}({\bf s})}{l_{m}(\g)}\;\;
\mbox{\rm for all $j \neq m$, $\;\;j,m \in \{ 1, \ldots, p\}$.}
$$
The open set $U(\g)$ is non-empty, since we assume that
$g.c.d.(l_j, l_{m})=1$  for $j \neq m$.
For ${\bf s} \in U(\g)$, we have
$$
{\rm Res}_{z = z_j} W({\bf s} + z \cdot \g) =
\frac{1}{(k_j-1)!}
\left( \frac{\partial}{\partial z} \right)^{k_j-1}
\frac{l_{j}({\bf s} + z \cdot \g)^{k_j}
P({\bf s} + z  \cdot \g)}{l_j^{k_j}
(\g) \cdot  \prod_{m =1}^p
l_{m}^{k_m}({\bf s} + z \cdot \g) }|_{z = z_j},
$$
where
$$
 z_j = - \frac{l_j({\bf s})}{l_j(\g)}.
$$

Let
$$
W({\bf s}) \cdot \prod_{j =1}^p l_j^{k_j}({\bf s})
= {P}({\bf s}) = {P}_0({\bf s}) + {P}_1({\bf s}),
$$
where ${P}_0({\bf s})$ is a uniquely determined homogeneous polynomial
and ${P}_0({\bf s})$ is an analytic at ${\bf 0}$ function
such that  $\mu({P}) = \mu({P}_0) < \mu({P}_1)$ and
$$
 \frac{P_0({\bf s})}{\prod_{j =1}^p l_j^{k_j}({\bf s})} =
C(W) \cdot {\cal  X}_{\L}({\bf s})
$$
(${\cal  X}_{\L}({\bf s})$ is the ${\cal X}$-function of the cone $\L =
\psi({\R}^r_{\geq 0})$).  We set
$$
R_0({\bf s}) : = \frac{P_0({\bf s})}{\prod_{j =1}^p {l}^{k_j}_j({\bf
s})}, \;\;
R_1({\bf s}) : = \frac{P_1({\bf s})}{\prod_{j =1}^p {l}_j^{k_j}({\bf
s})}.
$$
Then $\mu(W)=  \mu (R_0) < \mu (R_1)$. Moreover,
$\mu(W) = - {\rm dim} V_{\R} = r - {\rm rk}\, \Gamma$.
Define
$$
\tilde{R}_0({\bf s}):= 2\pi i \cdot
 \sum_{j=1}^{p_1} {\rm Res}_{z = z_j}
R_0({\bf s}+ z\cdot \g)
$$
and
$$
\tilde{R}_1({\bf s}):=
2\pi i \cdot  \sum_{j=1}^{p_1} {\rm Res}_{z = z_j}
R_1({\bf s}+ z\cdot \g).
$$
By Proposition \ref{mult4}, we
have
$\mu (\tilde{R}_1) \geq 1+ \mu(R_1)  \geq 2 + \mu(R_0)= 1+ {\rm
rk}\, \tilde{\Gamma}  - r $.

We claim
$$
\tilde{R}_0({\bf s}) = 2 \pi i \cdot C(W) {\cal  X}_{\tilde{\L}}
(\tilde{\psi}({\bf s}))
$$
in particular  $\mu (\tilde{R}_0) = \mu(R_0) + 1 = {\rm
rk}\, \tilde{\Gamma}$. Indeed, repeating for ${\cal  X}_{\L}(\psi({\bf s}))$
the same arguments as for $W({\bf s})$ we obtain
$$
\int_{{\rm Re}(z) = 0} {\cal  X}_{\L}(\psi({\bf s} + z \cdot \g)) dz
- \int_{{\rm Re}(z) = \d} {\cal  X}_{\L}(\psi({\bf s} + z \cdot \g)) dz
$$
$$
= 2\pi i \cdot \sum_{j=1}^{k_1} {\rm Res}_{z = z_j}
{\cal  X}_{\L}(\psi({\bf s} + z_j \cdot \g)).
$$
Moving the contour of integration
${\rm Re}(z) = \d$, by residue theorem,
we obtain
$$
\int_{{\rm Re}(z) = \d} {\cal  X}_{\L}(\psi({\bf s} + z \cdot \g)) dz =0.
$$
On the other hand,
$$
{\cal  X}_{\tilde{\L}}(\tilde{\psi}({\bf s})) =
\frac{1}{2\pi i}\int_{{\rm Re}(z) = 0}
{\cal  X}_{\L}(\psi({\bf s} + z \cdot \g)) dz
$$
(see Theorem \ref{char1}).

By \ref{mult2}(iii), using the decomposition
$$
\tilde{W}({\bf s}) = \tilde{W}_{\d}({\bf s}) + \tilde{R}_0({\bf
s}) + \tilde{R}_1({\bf s})
$$
and our assumption $\mu(\tilde{W}_{\d}) \geq
1 + {\rm rk}\, \tilde{\Gamma} - r$,
we obtain that $\mu (\tilde{W}) = \mu (\tilde{R}_0) = {\rm
rk}\, \tilde{\Gamma}  - r$.

By \ref{desc2}, the linear forms
\[ h_{m,j}({\bf s}):= l_{m}({\bf s} + z_j \cdot \g)
= l_{m}({\bf s}) -
\frac{l_j({\bf s})}{l_{j}(\g)} l_{m}(\g) \]
and the analytic in the domain $U(\g)$ functions
\[ {\rm Res}_{z = z_j} W({\bf s} + z \cdot \g), \;\;
{\rm Res}_{z = z_j} R_0({\bf s} + z \cdot \g)  \]
descend to $\tilde{V}_{\C}$.
For any $j \in \{ 1, \ldots, p_1\}$, let us denote
$$
Q_j ({\bf s}) = \prod_{m \neq j, m=1}^p h_{m,j}^{k_m}({\bf s}).
$$
It is clear that
$$
Q_j^{k_j}({\bf s}) \cdot {\rm Res}_{z = z_j} W({\bf s} + z \cdot \g)\;
\;\mbox{\rm and}\;\;
Q_j^{k_j}({\bf s}) \cdot {\rm Res}_{z = z_j} R_0({\bf s} + z \cdot \g)
$$
are analytic at ${\bf 0}$ and $\Phi({\bf s})$ divides
each $Q_j ({\bf s})$. So we obtain that
$$
 \tilde{W}({\bf s}) \prod_{j=1}^{p_1} Q_j^{k_j}({\bf s})
  =  \left( \tilde{W}_{\d}({\bf s}) +
  2 \pi i\cdot
\sum_{j=1}^{p_1} {\rm Res}_{z = z_j} W({\bf s} + z \cdot \g) \right)
\prod_{j=1}^{p_1} Q_j^{k_j}({\bf s})
$$
and
$$
 \tilde{R}_0({\bf s}) \prod_{j=1}^{p_1} Q_j^{k_j}({\bf s})
  =  \left(   2 \pi i\cdot
\sum_{j=1}^{p_1} {\rm Res}_{z = z_j} R_0({\bf s} + z \cdot \g) \right)
\prod_{j=1}^{p_1} Q_j^{k_j}({\bf s})
$$
are  analytic at ${\bf 0}$.

Define the set $\{ \tilde{l}_1({\bf s}), \ldots, \tilde{l}_q({\bf s}) \}$ as
a subset of pairwise coprime elements
in the set of homogeneous linear forms $\{ h_{m,j}({\bf s}) \}$ $(m
\in \{1, \ldots, p\}, \; j \in \{1, \ldots, p_1\})$  such that
there exist positive integers $n_1, \ldots, n_q$ and a
representation of the meromorphic functions $\tilde{W}({\bf s})$
and  $\tilde{R}_0({\bf s})$
as quotients
\[ \tilde{W}({\bf s}) =
\frac{\tilde{P}({\bf s})}{\prod_{j =1}^q \tilde{l}^{n_j}_j({\bf s})},\;\;
\tilde{R}_0({\bf s}) =
\frac{\tilde{P}_0({\bf s})}{\prod_{j =1}^q \tilde{l}^{n_j}_j({\bf s})},\]
where $\tilde{P}({\bf s})$ is analytic at ${\bf 0}$,
$\tilde{P}_0({\bf s})$ is a homogeneous polynomial  and
none of the forms
$\tilde{l}_1({\bf s}), \ldots,
\tilde{l}_q({\bf s})$ vanishes for $ {\bf s} \in {\R}_{>0}^r$
(the last property can be achieved, because both functions
$\tilde{W}({\bf s})$ and $\tilde{R}_0({\bf s})$ are  analytic in $U$).

Define
$$
\tilde{P}_1({\bf s}) = \left( \tilde{W}_{\d}({\bf s}) +
\tilde{R}_1({\bf s}) \right)
\cdot \prod_{j =1}^q \tilde{l}^{n_j}_j({\bf s}).
$$
Then
$$
\tilde{P}({\bf s}) = \tilde{P}_0({\bf s}) +
\tilde{P}_1({\bf s})
$$
where $\tilde{P}_0({\bf s})$ is a homogeneous polynomial
and $\tilde{P}_1({\bf s})$ is an analytic at ${\bf 0}$ function
such that  $\mu(\tilde{P}) = \mu(\tilde{P}_0) < \mu(\tilde{P}_1)$ and
$$
 \frac{\tilde{P}_0({\bf s})}{\prod_{j =1}^q
\tilde{l}_j^{n_j}({\bf s})} =
2\pi i \cdot C(W) \cdot
{\cal  X}_{\tilde{\L}}(\tilde{\psi}({\bf s})).
$$

\section{Main theorem}

Let us set
$$W_{\S}({\bf s}) :=
Z_{\S}(\varphi_{\bf s} +\p_{\S}) = Z_{\S}(s_1 +1, \ldots, s_r +1).$$
By Theorem \ref{convergence}, $W_{\S}({\bf s})$ is
an  analytic function in the
domain ${\rm Re}({\bf s}) \in {\bf R}^r_{>0}$.

\begin{theo}
The analytic function $W_{\S}({\bf s})$
is good with respect to the lattice
$M^G \subset PL(\S)^G = {\bf Z}^r$.
\label{analytic.cont}
\end{theo}

\noindent
{\em Proof.}
By Theorem \ref{poiss}, we have the following integral
representation for $Z_{\S}({\bf s}) $ in the domain
${\rm Re}({\bf s}) \in {\bf R}^r_{>1}$
$$
Z_{\Sigma}({\bf s})=\frac{1}{(2\pi )^t b_S(T)}
\int_{(T({\bf A}_K)/T(K))^*}\hat{H}_{\S}(\chi, -{\bf s})d\chi
$$
We need only to consider characters $\chi$
which are trivial
on the maximal compact subgroup
${\bf K}_T\subset T^1({\bf A}_K)$, because
for all other characters the Fourier transform
$\hat{H}_{\S}(\chi, -{\bf s})$ vanishes.
Choosing a non-canonical splitting of characters
corresponding to some
splitting of the sequence
$$
0\ra T^1({\bf A}_K)\ra T({\bf A}_K)\ra T({\bf A}_K)/T^1({\bf A}_K)\ra 0
$$
we obtain
$$
Z_{\Sigma}({\bf s})=\frac{1}{(2\pi )^t b_S(T)}
\int_{(T({\bf A}_K)/T^1({\bf A}_K))^*}d\chi_y
\int_{(T^1({\bf A}_K)/T(K))^*}\hat{H}_{\S}(\chi, -{\bf s})d\chi_l
$$
We have an isomorphism
$M^G_{\R}\simeq(T({\bf A}_K)/T^1({\bf A}_K))^*$ and the
measure $d\chi_y $ coincides with the usual Lebesgue measure
on $M^G_{\R}$.
Recall that a character $\chi\in (T({\bf A}_K)/T(K))^*$
defines $r$ Hecke characters $\chi_1,...,\chi_r$
of the idele groups ${\bf G}_m({\bf A}_{K_j})$.
In particular, we get $r$ characters $\chi_{1,y},...,\chi_{r,y}$.
We have an embedding $M^G\subset PL(\S)^G$,
which together with explicit formulas for Fourier transforms
of local heights shows that the integral
$$
A_{\S}({\bf s},\chi_y):=\frac{1}{b_S(T)}
\int_{(T^1({\bf A}_K)/T(K))^*}
\hat{H}_{\S}(\chi, -({\bf s}+{\bf 1}))d\chi_l
$$
is a function on $PL(\S)^G_{\C}$ and we have
$$
A_{\S}({\bf s},\chi_y)=A_{\S}({\bf s}+i{\bf y})= A_{\S}(s_1+iy_1,...,s_r+iy_r).
$$

Denote by $\Gamma:=M^G$ the lattice of $K$-rational characters
of $T$. Let $t$ be the rank of $\Gamma$.
The case $t =0$ corresponds to
an anisotropic torus $T$. It has been considered already in
\cite{BaTschi}. So we assume $t >0$.

For any element $\g \in \Gamma\subset {\Z}^r$ we
denote by $l(\g)$ the number of its
coordinates which are zero $(0 \leq l(\g) \leq r)$.
Let $l(\Gamma)$ be the minimum of $l(\g)$ among $\g \in \Gamma$.
 Notice that
$l(\Gamma) \leq r - t - 1$. Indeed, if we had
$l(\Gamma) \geq r -t$, then $M^G$ would be contained in
the intersection of $r - t$ linear coordinate hyperplanes
$s_j = 0$ (the latter contradicts the condition $M^G_{\R} \cap
{\bf R}^r_{\geq 0} = 0$).
We can always choose a
${\Z}$-basis $\g^1, \ldots ,\g^t$ of $\Gamma$ in such a
way that $l(\Gamma) = l(\g^u)$ $(u =1, \ldots, t)$. Without loss
of generality we assume
that $\Gamma$ is contained in the intersection of  coordinate
hyperplanes $s_j = 0$, $j \in \{1, \ldots, l(\Gamma) \}$.
We set
$$\Phi({\bf s}) := \prod_{j =1}^{l(\Gamma)} s_j. $$
For any $u \leq t$ we define a subgroup $\Gamma^{(u)} \subset
\Gamma$ of rank $u$ as
$$\Gamma^{(u)}:= \bigoplus_{j =1}^u {\Z}<\g_j>.$$
We introduce
some auxiliary functions
$$
W^{(u)}_{\S}({\bf s}) =
\int_{\Gamma^{(u)}_{\R}}
A_{\S}({\bf s}+i{\bf y}^{(u)}){\bf dy}^{(u)}
$$
where $ {\bf dy}^{(u)}$ is the induced measure
on $\Gamma^{(u)}_{\R}\subset PL(\S)^G_{\R}$.
Denote $V^{(u)}_{\C} = {\C}^r/\Gamma_{\C}^{(u)}$.
We prove by induction
that $W^{(u)}_{\S}({\bf s})$ is good with respect to
$\Gamma^{(u)} \subset {\Z}^r$.

By \ref{infconver}, $W^{(u)}_{\S}({\bf s})$ is an analytic function
in the domain ${\rm Re}({\bf s}) \in {\R}^r_{>0}$.

There exist $\d_1,...,\d_t > 0$ such that the
integral
$$
 \int_{{\rm Re}(z) = \d_u}
\Phi({\bf s}) \cdot W_{\S}^{(u-1)}({\bf s} + z \cdot \g^u) dz
$$
converges absolutely and uniformly
in an  open neighborhood of ${\bf 0}$.
This can be seen as follows:
For any $\e$ with  $0<\e<1/rd'$,
where $d'=\dim M_{\R,\infty}$, we can
choose a
ball $B_{e_1}\subset \R$ of radius $e_1$ around ${\bf 0}$ such that
for any ball
$B_{e_2}\subset B_{e_1} $ of radius $e_2$ ($0<e_2 <e_1$)
around  $0$ we can uniformly bound the Hecke $L$-functions
$L_{K_j}(s_j+1,\chi_j)$
appearing in $\hat{H}_{\S}(\chi,{\bf s})$  by
$$
c_j(e_2 )(\|y(\chi_j)\|+ |{\rm Im}(s_j)|+1)^{\e}
$$
with some constants $c_j(e_2 )$ for all
${\bf s}$ in the domain ${\rm Re}(s_j)\in B_{e_1}\backslash B_{e_2}$
for $j=1,...,r$ (see \ref{estim}).
By \ref{infconver}, this assures the
absolute and uniform convergence
of the integral
$$
\int_{\Gamma^{(u)}_{\R}}
A_{\S}({\bf s}+i{\bf y}^{(u)}){\bf dy}^{(u)}
$$
for all ${\bf s}$ contained
in a compact in $\C^r$ such that
${\rm Re}(s_j)\in B_{\e_1}\backslash B_{\e_2}$
for $j=1,...r$.
We know that the  coordinates $\g_j^u$ of
the vectors $\g^u=(\g_1^u,...,\g^u_r)\in \Z^r$
are not equal to zero for $l(\Gamma )<j \le r$.
Therefore, we  can now choose such real
$\d_u>0$ that  $\d_u\g^u_j$
are all contained in the {\em open} ball $B_{e_1}$.
So there must exist some $e_2>0$ such that
$\d_u\g^u_j\not\in B_{e_2}$ for all $u=1,...,t$ and all
$l(\Gamma )<j\le r$. It follows that there
exists an open neighborhood of
${\bf 0}$, such that for all
${\bf s}$ contained in this neighborhood
we have ${\rm Re}(s_j+\d_u\g^u_j)\in
B_{e_1}\backslash B_{e_2}$ for all
$l(\Gamma )<j\le r$. Since we remove the remaining poles
by multiplying with
$\Phi({\bf s})$ we obtain  the absolute
and uniform convergence of $W^{(u)}_{\S}({\bf s})$ to a
holomorphic function in ${\bf s}$ in this neighborhood.

Moreover, the
multiplicity of the meromorphic function
$$
\tilde{W}_{\d_u}^{(u)}({\bf s}): =   \int_{{\rm Re}(z) = \d_u}
 W^{(u-1)}_{\S}({\bf s} + z \cdot \g_u) dz
$$
at  ${\bf 0}$ is at least $1 + {\rm rk}\, {\Gamma} - r \geq 1 +
{\rm rk}\, {\Gamma}^{(u)} - r$.

We apply Theorem \ref{desc3}. The property (iii) follows from
estimates \ref{m.estim} and \ref{infconver}.
This concludes
the proof.

\hfill $\Box$

\begin{theo} Denote by  $\hat{H}_{\S,S}(\chi,-{\bf s})$
the multiplicative Fourier transform of the height function
with respect to the measure $\omega_{\Omega, S}$ (see \ref{can.meas}).
The principal coefficient $C(\S)$
of
$$
A_{\S}({\bf s})=\frac{1}{b_S(T)} \int_{(T^1({\bf A}_K/T(K){\bf K}_T)^*}
\hat{H}_{\S,S}(\chi_l,-{\bf s})d\chi_l
$$
at $s_1=...=s_r=1$ is equal to
$\beta({\bf P}_{\S}) \tau_{\cal K}({\bf P}_{\S})$.
\label{beta.tau}
\end{theo}

{\em Proof.} We follow closely the exposition
of the proof of theorem 3.4.6 in \cite{BaTschi}.
Since $M^G\hookrightarrow PL(\S)^G$ we have an embedding
of characters
$$
(T({\bf A}_K)/T^1({\bf A}_K))^* = M^G_{\R}\hookrightarrow
\prod_{j=1}^r ({\bf G}_m({\bf A}_{K_j})/{\bf G}_m^1({\bf A}_{K_j}))^*.
$$
Recall that the kernel of
$$
a^*\,:\, (T({\bf A}_K)/T({\bf A}_K))^* \ra
\prod_{j=1}^r ({\bf G}_m({\bf A}_{K_j})/{\bf G}_m({\bf A}_{K_j}))^*
$$
is dual to the obstruction group to weak approximation
$A(T)=T({\bf A}_K)/\overline{T(K)}$.
We have a splitting
$$
 \overline{T(K)} =   \overline{T(K)}_S \times T(A_{K,S}).
$$
Here we denoted by
$\overline{T(K)}_S$  the image of $\overline{T(K)}$ in
$\prod_{v \in S} T(K_v)$ and
$T(A_{K,S})=T({\bf A}_K)\cap \prod_{v\not\in S}T(K_v)$.
The pole of the highest order
$r$ of $\hat{H}_{\S,S}(\chi_l,-{\bf s}) $ at $s_1=...=s_r=1$ appears
from characters $\chi_l$ such that the corresponding
$\chi_1,...,\chi_r$
are trivial characters of the groups
${\bf G}_m({\bf A}_{K_j})/{\bf G}_m({K_j})$,
i.e., $\chi_l$ is a character of the finite group
$A(T)=\prod_{v\in S}T(K_v)/\overline{T(K)}_S$,
and is trivial on the group
$T({\bf A}_{K,S})$.

For ${\bf s}\in \R_{>1}^r$ we can again apply the Poisson
formula to $A(T)$. By \ref{weak1}, the order of $A(T)$ equals
$\beta({\bf P}_{\S})/i(T)$. We obtain
$$
\frac{1}{b_S(T)}
\sum_{\chi \in (A(T))^* }
\hat{H}_{\S,S}(\chi_l, - {\bf s})=
 \frac{\beta({\bf P}_{\Sigma})}{i(T) b_S(T)}
\int_{\overline{T(K)}} H_{\S}(x,-{\bf s})
\omega_{\Omega,S} \]
(see  \ref{weak1}).
We restrict to the line $s_1=...=s_r=s$ and we want to compute the
limit
$$
\lim_{s\ra 1} (s-1)^r \int_{\overline{T(K)}} H_{\S}(x,-{\bf s})
\omega_{\Omega,S}.
$$
We have
\begin{equation}
 \int_{\overline{T(K)}} H_{\S}(x,-{\bf s})
\omega_{\Omega,S} =
\label{const1}
\end{equation}
$$
= \int_{\overline{T(K)}_S}
\prod_{v \in S} H_{\S,v}(x_v,-{\bf s}) \omega_{\Omega,v}
\cdot
\prod_{v \not\in S} \int_{T(K_v)} H_{\S,v}(x_v,-{\bf s}) d\mu_v
$$
(recall that   $\omega_{\Omega,v}=\prod_{v\in {\rm Val}(K)} d\mu_v$
and
$d\mu_v = L_v(1,T;E/K) \omega_{\Omega,v}$ for all
$v$ and $L_v(1,T;E/K) =1$ for $v \in S$).

{}From our calculations of the Fourier transform of local
height functions for $v \not\in S$ (\ref{loc-int}),
we have
\begin{equation}
\prod_{v \not\in S} \int_{T(K_v)}
H_{\S,v}(x_v,-{\bf s}) d\mu_v
=
\label{const2}
\end{equation}
\[ =  L_S(s, T;E/K)
\cdot L_S(s, T_{NS}; E/K) \prod_{v \not\in S}
Q_{\S}(q_v^{-s}, \ldots, q_v^{-s}). \]
By \ref{p-function},
\[ \prod_{v \not\in S}
Q_{\S}(q_v^{-s}, \ldots, q_v^{-s}) \]
is an absolutely convergent Euler product for $s =1$.
Moreover, the limits
$$
\lim_{s\ra 1} (s-1)^t L_S(s, T;E/K)
$$
$$
\lim_{s\ra 1} (s-1)^{(r-t)} L_S(s, T_{NS};E/K)
$$
exist and equal the non-zero constants
$l_S(T)$ and $l^{-1}_S({\bf P}_{\S})$ ($r=t+k$).
By \ref{badreduction},
\[   \int_{\overline{T(K)}_S}
\prod_{v \in S} H_{\S,v}(x_v,-{\bf s}) \omega_{\Omega,v} \]
is  absolutely convergent for $s_1=,,,=s_r=1$.
Using (\ref{const1}) and (\ref{const2}), we obtain:
\begin{equation}
 \lim_{s \ra 1} (s-1)^r
\int_{\overline{T(K)}} H_{\S}(x,-{\bf s})
\omega_{\Omega,S} =
\label{const3}
\end{equation}
\[ = \frac{l_S(T)}{l_S({\bf P}_{\S})}
\int_{\overline{T(K)}_S}
\prod_{v \in S} H_{\S,v}(x_v,-{\bf s}) \omega_{\Omega,v} \cdot
 \prod_{v \not\in S}
Q_{\S}(q_v^{-1}, \ldots, q_v^{-1}).\]
Now recall (\ref{loc-int}),
that for $v\not\in S$ we have
$$
Q_{\S}(q_v^{-1}, \ldots, q_v^{-1})=
\int_{T(K_v)}L_v^{-1}(1,T_{NS};E/K)
H_{\S,v}(x_v,-{\bf 1})\omega_{\Omega,v}.
$$
It was
proved in \cite{BaTschi} Proposition 3.4.4
that the restriction of the $v$-adic
measure $\omega_{{\cal K},v}$ to $T(K_v) \subset
{\bf P}_{\Sigma}(K_v)$ coincides with the measure
\[ H_{\S,v}(x, -{\bf 1}) \omega_{\Omega,v}. \]
Here  ${\cal K}$ is the
canonical sheaf on the toric variety ${\bf P}_{\Sigma}$ metrized
as above.

We also have
\begin{equation}
 \int_{\overline{T(K)}_S}
\prod_{v \in S} H_{\S,v}(x_v,-{\bf 1}) \omega_{\Omega,v} =
\int_{\overline{T(K)}_S}
\prod_{v \in S} \omega_{{\cal K},v}.
\label{const6}
\end{equation}

Using the splitting
$\overline{T(K)} =  \overline{T(K)}_S \times T(A_{K,S})$
and multiplying the above equations we get
$$
\int_{\overline{T(K)}} \omega_{{\cal
K},S} =  \int_{\overline{T(K)}_S} \prod_{v \in S} \omega_{{\cal K},v}
\cdot  \prod_{v \not\in S} \int_{{T(K_v)}}
L_v^{-1}(1,T_{NS};E/K) \omega_{{\cal K},v}.
$$
On the other hand, it was proved in \cite{BaTschi} Proposition
3.4.5 that we have
\[  \int_{\overline{T(K)}} \omega_{{\cal
K},S} = \int_{\overline{{\bf P}_{\S}(K)}}
\omega_{{\cal K},S} = b_S({\bf P}_{\S}). \]

Therefore,
\[ b_S({\bf P}_{\S}) =
\int_{\overline{T(K)}_S}
\prod_{v \in S} H_{\S}(x,-\varphi_{\Sigma}) \omega_{\Omega,v}
\cdot \prod_{v \not\in S}
Q_{\S}(q_v^{-1}, \ldots, q_v^{-1})
.\]

Collecting terms, we obtain
$$
C(\S )= \frac{\beta({\bf P}_{\S})}{i(T)  b_S(T)} \cdot
\frac{l_S(T)}{l_S({\bf P}_{\S})} \cdot b_S({\bf P}_{\S}).
$$
By  \ref{tamagawa1} and \ref{tamagawa}, we have
the following equality
$$
 i(T) b_S(T) = h(T) l_S(T).
$$
It remains to notice that
we have an exact sequence of lattices
$$
0\ra M^G\ra PL(\S)^G\ra {\rm Pic}({\bf P}_{\S})\ra H^1(G,M)\ra 0
$$
and that the number $h(T)= |H^1(G,M)|$
appears in the integral formula
for the ${\cal X}$-function of the
cone ${\L }_{\rm eff}\subset {\rm Pic}({\bf P}_{\S})$.
We apply Theorem \ref{char0} and obtain that
$$
W_{\S}({\bf s})=\frac{1}{(2\pi)^t b_S(T)}
\int_{M^G_{\R}}A_{\S}({\bf s}+i{\bf y}){\bf dy}
$$
is good  with respect to the lattice $M^G$ and that
$$
C(\S)= \beta( {\bf P}_{\S})\tau_{\cal K}({\bf P}_{\S})
$$
is the principal coefficient of $W_{\S}({\bf s})$ at ${\bf 0}$.
\hfill $\Box $

\begin{theo} There exists a $\delta >0$ such that the
height zeta-function $\zeta_{\S}(s)$ obtained by restiction of
the zeta-function
$Z_{\S}({\bf s})$ to the complex line $s_j = \varphi (e_j)=s$ for
all $j=1,...,r$ has a representation of the form
$$
\zeta_{\S}(s)= \frac{\Theta(\S)}{(s-1)^k} + \frac{g(s)}{(s-1)^{k-1}}
$$
with $k= r-t = {\rm rk}\, {\rm Pic}({\bf P}_{\S})$ and  some
holomorphic function $g(s)$ in the domain
${\rm Re}(s)>1-\delta$. Moreover, we have
$$
\Theta(\Sigma)  = \alpha({\bf P}_{\Sigma})\beta({\bf P}_{\Sigma})
\tau_{\cal K}({\bf P}_{\Sigma}).
$$
\end{theo}

{\it Proof.} Since $W_{\S}({\bf s})$ is good with respect to the
lattice $M^G \subset {\Z}^r$, we have the following representation
of $W_{\S}({\bf s})$ in a small open neighborhood of ${\bf 0}$:
$$
W_{\S}({\bf s}) = \frac{ P({\bf s})}{\prod_{j =1}^p l_j^{k_j}({\bf s})}
$$
where $P({\bf s}) = P_0({\bf s}) + P_1({\bf s})$, $\mu(P_1) >
\mu(P_0)$ and
$$
\frac{P_0({\bf s})}{\prod_{j =1}^p l_j^{k_j}({\bf s})} =
\beta({\bf P}_{\S}) \tau_{\cal K}({\bf P}_{\S})
\cdot {\cal X}_{\L_{\rm eff}}(\psi({\bf s})),
$$
where ${\cal X}_{\L_{\rm eff}}$ is the ${\cal X}$-function of the cone
$\L_{\rm eff} = \psi({\R}^r_{\geq 0}) \subset {\rm Pic}({\bf P}_{\S})_{\R}$.

If we restrict
$$
\frac{P_0({\bf s})}{\prod_{j =1}^p l_j^{k_j}({\bf s})}
$$
to the line
$s_j = s - 1$ $(j =1, \ldots, r)$, then
we get the meromorphic function $\Theta (s-1)^{-k}$ with
$\Theta = \alpha({\bf P}_{\Sigma})\beta({\bf P}_{\Sigma})
\tau_{\cal K}({\bf P}_{\Sigma})$.
Moreover, the order of the pole at $s =1$
of the restriction of
$$
\frac{P_1({\bf s})}{\prod_{j =1}^p l_j^{k_j}({\bf s})}
$$
to the line
$s_j = s - 1$ $(j =1, \ldots, r)$ is less than $k$.
Therefore, this restriction
can be written as $g(s)/(s-1)^{k-1}$
for some analytic at $s =1 $ function
$g(s)$.

\begin{coro}
{\rm
Let  $T$ be an algebraic
torus and ${\bf P}_{\Sigma}$ its smooth
projective compactification.
Let $k$ be
the rank of ${\rm Pic}({\bf P}_{\Sigma})$.
Then the number
of $K$-rational points  $x \in T(K)$ having the anticanonical
height $H_{{\cal K}^{-1}}(x) \leq B$ has the asymptotic
\[ N(T,{\cal K}^{-1}, B) = \frac{\Theta(\Sigma)}{(k-1)!}
\cdot B (\log B)^{k-1}(1+o(1)), \hskip 0,3cm B\ra \infty.\]
}
\end{coro}

\noindent
{\em Proof.}
We apply a Tauberian theorem to $\zeta_{\S}(s)$.
\hfill $\Box$

\end{document}